\documentclass[10pt]{amsart}

\usepackage{nomencl}

\usepackage{mathrsfs,bbm}
\usepackage{epic,eepic,graphics}
\usepackage{pstricks}
\usepackage{pst-plot}
\usepackage[latin1]{inputenc}
\usepackage{verbatim}

\def\varnothing{\emptyset}

\newcommand{\pr}[1]{{\rm {\mathbf P}}\left(#1\right)}

\newcommand{\expect}[1]{{\rm {\mathbf E}}\left[#1\right]}

\newcommand{\beeq}{\begin{equation}}
\newcommand{\eneq}{\end{equation}}
\newcommand{\bear}{\begin{eqnarray}}
\newcommand{\enar}{\end{eqnarray}}
\newcommand{\bearno}{\begin{eqnarray*}}
\newcommand{\enarno}{\end{eqnarray*}}

\newtheorem{thm}{Theorem}[section]
\newtheorem{prop}[thm]{Proposition}
\newtheorem{property}[thm]{Property}
\newtheorem{cor}[thm]{Corollary}
\newtheorem{lem}[thm]{Lemma}
\newtheorem{defn}[thm]{Definition}

\newcommand{\cp}{\cite}

\newcommand{\paren}[1]{\left(#1\right)}
\newcommand{\dw}[1]{\mathbf{d}\left[#1\right]}
\newcommand{\ww}[1]{\mathbf{w}_T\left(#1\right)}

\newcommand{\ar}{\rightarrow}
\newcommand{\Ar}{\Rightarrow}

\newcommand{\wk}{\overset{\w}{\longrightarrow}}

\newcommand{\lar}{\longrightarrow}

\DeclareMathOperator{\w}{{\mathbf w}}

\newcommand{\BP}{{\mathbf P}}
\newcommand{\BE}{{\mathbf E}}

\newcommand{\BPr}{{\bf P}^r}

\newcommand{\f}{\frac}

\newcommand{\qq}{\qquad}
\newcommand{\cd}{(\cdot )}

\newcommand{\BN}{{\mathbb N}}

\newcommand{\BR}{{\mathbb R}}
\newcommand{\Rp}{{\mathbb{R}}_+}

\newcommand{\DM}{\mathbf{D}\left([0,\infty),\bM\right)}

\newcommand{\R}{\mathcal{R}}

\newcommand{\M}{\mathbf{M_1}}

\newcommand{\bM}{\mathbf{M}}

\newcommand{\K}{\mathbf{K}}

\newcommand{\Er}{\overline{E}^r}

\newcommand{\Lbr}{\overline{\mathcal{L}}^r}

\newcommand{\vbr}{\overline{B}^r}

\newcommand{\lbr}{\overline{D}^r}

\newcommand{\Sr}{\overline{S}^r (t,t+h)}
\newcommand{\Sbr}{\overline{S}^r}

\newcommand{\Z}{\mathcal{Z}}

\newcommand{\Zbr}{\overline{\mathcal{Z}}^r}

\def\etal{{\em et al. }}

\title[Processor Sharing Queues with Impatience]{Fluid Limits for
Processor Sharing Queues with~Impatience}

\author[H.C. Gromoll]{H.\ Christian
Gromoll${}^*$}\thanks{${}^*$Research supported in part by an NSF
Mathematical Sciences Postdoctoral Research Fellowship, a European Union
Marie Curie Postdoctoral Research Fellowship, and EURANDOM}
\address[Christian Gromoll]{Department of Mathematics,
Stanford University  450 Serra Mall, Stanford, CA 94305-2125, USA}
\email{gromoll@math.stanford.edu}

\author[Ph. Robert]{Philippe Robert}
\address[Ph. Robert]{INRIA{-}Rocquencourt, RAP project, Domaine de
Voluceau, BP 105, 78153
 Le~Chesnay, France}
\email{Philippe.Robert@inria.fr}
\urladdr{http://www-rocq.inria.fr/$\tilde{\ }$robert}

\author[B. Zwart]{Bert Zwart${}^\dag$}\thanks{${}^\dag$Research supported by an NWO-VENI grant}
\address[Bert Zwart]{Eindhoven University of Technology,
Department of Mathematics and Computer Science,
HG 9.35,
P.O. Box 513,
5600 MB Eindhoven,
the Netherlands}
\email{zwart@win.tue.nl}


\begin{document}

\begin{abstract}
We investigate a  processor sharing queue with renewal  arrivals and generally distributed
service  times.   Impatient jobs  may  abandon the  queue,  or  renege, before  completing
service. The random time representing a  job's patience has a general distribution and may
be dependent on its initial service  time requirement. A scaling procedure that gives rise
to a  fluid model with nontrivial yet  tractable steady state behavior  is presented. This
fluid model model captures many essential features of the underlying stochastic model, and
it is used to  analyze the impact of impatience in processor  sharing queues.  
\end{abstract}

\maketitle

\hrule

\tableofcontents

\vspace{-8mm}

\hrule

\bigskip

\section{Introduction}\label{s.intro}
\subsection*{Processor-Sharing Policy and Impatience}
Processor Sharing  (PS) policies  were originally  proposed as models  of time  sharing in
computer operating systems.   Recently, generalizations of this discipline  have been used
to  describe data  transfers in  congested routes  through the  Internet, see  Roberts and
Massoulié~\cite{RoMa2000}  and  Kelly  and  Williams~\cite{KeWi2004}  and  the  references
therein. This has created considerable renewed interest in the analysis of PS policies.

\nocite{BonMas2001,Bra2005,deVKonLee2001,GrWi2005,KeWi2004,LaBeSr2004,Massoulie,Mas2005,MoWal2000,RoMa2000}
%Too long list of papers. Keep them but do not list the numbers (PhR).

This paper studies  the behavior of a $GI/GI/1$ queue serving  impatient jobs according to
the PS policy: if there are $N$  jobs in the queue, each job receives simultaneous service
at rate $1/N$.  An {\em impatient job} has a random {\em initial lead time} in addition to
its service  time. Such  a job has  a {\em deadline}  equal to  its arrival time  plus its
initial lead  time; if the  job has  not completed service  when the deadline  expires, it
abandons  the queue  (or {\em  reneges})  and therefore  does not  complete service.   For
example,  the timeout  of  a TCP  flow  through the  Internet  can be  thought  of as  the
expiration of a random deadline and subsequent reneging of the flow.

The  impact of  impatience on  PS  queues is  larger than  for  First In  First Out  (FIFO)
queues.  A typical  job that  abandons a  FIFO queue  will do  so while  waiting  to begin
service. In contrast,  a job that abandons  a PS queue will have  already received partial
service.  Since this  partial service  is wasted, impatience  may  create a
significant overhead for a PS server. 

There is a large literature on  queueing models with impatience under the FIFO discipline.
An  early  paper  by  Barrer  \cite{B57}  considers  an  example  arising  in  a  military
application.  Stanford~\cite{Stanford} is  a survey of the literature  in this domain (see
also  Stanford~\cite{S79} and  Boots and  Tijms~\cite{BT98}).  This  body of  work focuses
primarily  on  exact performance  analysis.   Ward  and  Glynn \cite{WG04}  have  recently
obtained  a diffusion  approximation for  single channel  queues. There  are  also various
studies of  multi-server queues with  abandonments, motivated by call  center applications;
see the survey by Gans {\em et al.} \cite{GKM02} and references therein.

There is some related literature treating other  policies, but in the context of {\em soft
deadlines}.  Jobs with  soft deadlines are not impatient; they remain  in the system until
completing service, even if their deadlines  have expired. In particular, these queues are
work conserving.   Results for such models describe  the extent to which  overdue jobs are
produced  by  the underlying  service  discipline,  without  the effect  of  abandonments.
Doytchinov       {\em       et       al.}~\cite{DoLeSh2001},      Kruk       {\em       et
al.}~\cite{KruLehShr2003,KruLehShr2004},   and    Yeung   and   Lehoczky~\cite{YeuLeh2004}
investigate  the heavy traffic  behavior of  various systems  using the  Earliest Deadline
First and  FIFO policies.  Gromoll and  Kruk \cite{GrKru2004} describes  the heavy traffic
behavior  of a  PS  queue incorporating  a  fairly general  structure  of soft  deadlines.
\nocite{KruLehShrYeu2001,KruLehShrYeu2004}
% Again, to avoid long lists (PhR).

For PS  queues with impatience  however, only  a few results  are known.  Coffman  {\em et
al.}~\cite{CPRW94} cover  the special  case of exponential  service times and  lead times,
where the lead time and service time are independent.  Guillemin {\em et al.}~\cite{GRZ04}
consider heavy tailed  service times, and obtain some results on  the reneging behavior of
large jobs  by analyzing the tail behavior  of the sojourn time  distribution.  Using some
approximations, Bonald and  Roberts \cite{BR03} analyze the steady state  of a system with
general service times and some dependence between service times and lead times.

\subsection*{Results of the Paper}
This paper analyzes the  PS queue with impatience by using fluid  limits.  The dynamics of
the system are represented as a measure valued process: the system state at time $t\geq 0$
is    represented     by    a    random    point     measure    ${{\mathcal{Z}}}(t)$    on
$(0,\infty]\times(0,\infty]$,  such   that  ${{\mathcal{Z}}}(t)$  has  a   point  mass  at
$(b,d)\in(0,\infty]\times(0,\infty]$ if and  only if there is a job in  the system at time
$t$  with residual  service  time $b$  and residual  lead  time $d$.   See Jean-Marie  and
Robert~\cite{JeanMarie:02}   and  Doytchinov   \etal~\cite{DoLeSh2001}  for   an  analogue
representation of  residual service times  in single server  queues. This setup  enables a
fairly general  analysis.  The case of a  general joint distribution of  service times and
initial lead  times, with possible dependence of  the two random variables  is included in
our setting.

Under mild assumptions, it  is shown that, with a convenient scaling,  a family of measure
valued  processes  associated  with  $({{\mathcal{Z}}}(t))$  is  tight  and  converges  in
distribution to some $(\zeta(t))$.  For $t\geq 0$, $\zeta(t)$ is a nonnegative measure on
$(0,\infty]\times(0,\infty]$, it  is the limit in  distribution of the  sequence of random
points  describing the queue.   This fluid  limit is  characterized as  the solution  of a
functional  Equation~\eqref{measurefluid}   which  can  be   viewed  as  a   time  changed
functional differential equation.

The  overloaded case $\rho>1$, which forms our main focus,  
presents a  nontrivial and  quite interesting  steady state
behavior. The  total fluid mass in  the system at  equilibrium (the fluid analogue  of the
total number  of jobs)  is shown to  be the  solution $z_\infty$ of  a simple  fixed point
equation~\eqref{fixedpointeq}.   Moreover,  the fluid  steady  state,  i.e.  the limit  of
$\zeta(t)$ as  $t$ goes  to infinity,  is a distribution  on $(0,\infty] \times (0,\infty]$  
which has  a simple expression~\eqref{Rully} in terms of $z_\infty$.

These results give  also a significant insight on the qualitative  properties of PS queues
with  impatience. An interpretation  of the  fixed point  equation~\eqref{fixedpointeq} is
given and  used to  analyze the total  number of jobs  in the  system and to  estimate the
fraction of jobs  that renege. The impact of  the variability of the service  times and of
the  lead  times  and  other  properties   of  this  queue  are  extensively investigated  
in  Gromoll \etal~\cite{Gromoll:02}.

In  contrast to  the models  studied  previously in  this domain,  the service  discipline
considered here  is {\em  not work conserving}.   For this  reason, analysis of  the fluid
model is more intricate. This is an  important difference from earlier work on standard PS
queues where the fact that the workload process coincides with that of FIFO discipline was
{\em a crucial ingredient} in the proof of the key results.  A different approach to prove
existence,  uniqueness,  and convergence  to  steady state  of  fluid  model solutions  is
proposed.  It is shown that there exists a {\em maximal} fluid model solution and by using
monotonicity arguments, the properties of the fluid limits can be investigated under quite
general assumptions.

\subsection*{Organization}
The paper is organized as follows. A detailed description of the model and the main results
is presented in Section 2. Qualitative properties of the fluid model are analyzed in Section 3. 
Section 4 is devoted to examples. Section 5 and 6 are concerned with convergence towards 
the fluid limit. Section 5 establishes tightness, and Section 6 characterizes limit points.

\section{Model Description and Results}\label{s.model}\setcounter{equation}{0}
This section gives  a detailed description of the stochastic  processes associated to this
queue and a summary of the main results.
%The fluid model analyzed here can be seen  as a  first order approximation, and our
%results concerning the properties of this fluid model and its validity as an
%approximation. 

\subsection{Stochastic model.}\label{s.stochModel}
The stochastic model consists of the following: a processor sharing server working at unit
rate  from an  infinite capacity  buffer, a  collection of  stochastic  primitives $E\cd$,
$\{B_i,D_i\}_{i=1}^\infty$  describing  respectively  the  process of  arrivals  and  the
services and the deadlines of the  customers, and a random initial condition specifying the
state of  the system at time  $0$. All random objects  are defined on  a probability space
$(\Omega,\mathscr{F},\BP)$ with expectation operator $\BE(\cdot)$.

%The time  evolution of  the system  is  described by  stochastic processes,
%defined in  terms of  the primitives and  initial condition  through a set  of descriptive
%equations.    

The {\em  exogenous arrival  process} $(E(t), t\geq  0)$ has  rate $\lambda{>}0$, it  is a
delayed renewal  process starting from  zero, with $i$th  jump time $U_i$. For  $t\geq 0$,
$E(t)$ is  the number of jobs  that arrive to the  buffer during $(0,t]$.   For $i\geq 1$,
$U_i$ is the  arrival time of job $i$; jobs  already in the buffer at  time $0$ are called
{\em initial jobs}.

For  $i\ge 1$,  the  {\em  service time}  $B_i$  is a  strictly  positive random  variable
representing the  amount of  processing time that  job $i$  requires from the  server. The
random variable  $D_i$ is  strictly positive and  determines the  deadline of job  $i$: it
represents the  maximum amount  of time that  job $i$  will stay in  the buffer.  Since it
arrives at time  $U_i$, its deadline is at  time $U_i+D_i$. It will abandon  the system at
this time if  it has not yet completed  service.  The random variable $D_i$  is called the
{\em initial lead time} of job $i$.

The model allows either the service time or the initial lead  time (but not both) to
be equal to infinity.  Therefore, the  random variable $(B_i,D_i)$ has values in the space
$\overline{\BR}_+^2=[0,\infty]\times[0,\infty]$.   Here,  $\overline{\BR}_+=[0,\infty]$ is
the usual  compactification of $\R_+$ with  the arithmetic extensions  $x+\infty = \infty$
for all $x\in\overline{\BR}_+$, $x\cdot\infty=\infty$ for $x>0$ and $0\cdot\infty=0$.  The
collection of Borel subsets of $\overline{\BR}_+^2$ is denoted by $\mathscr{B}$.  
Throughout the paper, it
is assumed that all sequences of services and deadlines $\{B_i,D_i\}_{i=1}^\infty$ are
independent  and  identically  distributed  (i.i.d.)   $\overline{\BR}_+^2$-valued  random
variables, and that their  common joint  distribution $\vartheta$  on $\overline{\BR}_+^2$
satisfies
\[
	\vartheta(\{0\}\times\overline{\BR})
			=\vartheta(\overline{\BR}\times\{0\}) =\vartheta(
			(\infty,\infty) ) =0.
%\label{a.vartheta} 
\]
Note that the random variables $B_i$ and $D_i$ may be dependent.  A generic random element
of  $\overline{\BR}_+^2$   with  distribution  $\vartheta$   will  be  denoted   $(B,D)$.   Let
$\rho=\lambda\BE[B]$  denote the {\em  traffic intensity}  of the  system.  It  is assumed
throughout that $\rho>1$,  that is, the server is nominally overloaded.   In this way, the
classical PS queue, the  infinite server queue, and mixtures of the  two are special cases
of this model.  (For example the $GI/GI/\infty$ queue corresponds to the case when service
times are equal to infinity.) 

% (Note  that $\rho=\infty$ for the $GI/GI/\infty$ model.)

\subsection*{Initial condition}  
The {\em  initial condition} specifies $Z(0)$, the  number of initial jobs  present in the
buffer at time zero, as well as the  service times and initial lead times of these initial
jobs.  Assume that  $Z(0)$ is a non-negative, integer valued  random variable.  The service
times and initial lead times for initial  jobs are the first $Z(0)$ elements of a sequence
$({B}_j^0,{D}_j^0)$     of     i.i.d.\     random     variables     taking    values     in
$\{(0,\infty]\times(0,\infty]\}  \setminus  (\infty,\infty)$  almost surely.   A  generic
random  element of  $\overline{\BR}_+^2$  distributed as  $({B}_0^0,{D}_0^0)$  will be  denoted
by $({B^0},{D^0})$.   Assume   that  the   expected  number  of   initial  jobs   is  finite:
$\BE[Z(0)]<\infty$.

\subsection*{Time Evolution of the Queue} 
%Several  stochastic  processes  describe how  the  system  evolves  due to  job  arrivals,
%processor sharing service,  and job impatience. 
For each  $t\ge 0$,  let $Z(t)$ denote  the number of  jobs in  the buffer (or  {\em queue
length}) at time $t$, and $S(t)$  denotes the {\em cumulative service time per job} provided
by the server up  to time $t$. Because of the processor sharing policy, the quantity
$S(t)$ is given by  
\begin{equation} \label{eq:seq-S}
S(t)  = \int_0^t \frac{1}{Z(s)}\,{\rm d}s,
\end{equation}
where the integrand is defined to be zero when the queue length equals zero.
If a job arrived at time $s\geq 0$ and is still present in the queue at $t\geq s$, at time
$t$ it has received the cumulative amount of processing time $S(t)-S(s)$.   

%Similarly, an initial job $j\le Z(0)$ still in the buffer  at time $t$  receives
%cumulative  service $S(t)$  by time  $t$.  According  to the processor sharing policy, a
%job present in  the buffer at time $t$ receives service at the instantaneous rate
%$1/Z(t)$.  The cumulative service  per job up to time $t$ can therefore be written 

Therefore, the {\em residual service time} at time $t$ of job $i\leq E(t)$
(and of initial job $j\leq Z(0)$) are given by
\[% \label{eq:seq-R}               
 B_i(t)  = (B_i - (S(t)-S(U_i)))^+,\qquad\text{and}\qquad
  B^0_j(t)  = (B^0_j - S(t))^+.
\]
Define the {\em lead time} at time $t$ of job $i\leq E(t)$ (and of
initial job $j\leq Z(0)$) by
\begin{equation}\label{li} 
D_i(t)  =(U_i+D_i-t)^+,\qquad\text{and}\qquad
D^0_j(t)  =(D^0_j-t)^+. 
\end{equation}
A  job's residual service  time is  the remaining  amount of  processing time  required to
fulfill its service  requirement; its lead time is the remaining  time until its deadline.
Job $i$ will  depart the system either  when its service requirement is  fulfilled or when
its deadline arrives, it will leave the system at time
\begin{equation}
\inf\{t\geq U_i: \min\{B_i(t),D_i(t)\}=0\}.
\nonumber
\end{equation}

The {\em state  descriptor} is a measure  valued process that keeps track  of the residual
service times and lead  times of all jobs in the buffer. For  job $i$, this information is
represented as a unit of mass  at the point $(B_i(t),D_i(t))\in\overline{\BR}_+^2$ at all times
$t\ge U_i$ such  that job $i$ is  still in the system.  Let  $\delta^+_{(x,y)}$ denote the
Dirac   point   measure   at   $(x,y)\in\overline{\BR}_+^2$   if   $\min\{x,y\}>0$,   otherwise
$\delta^+_{(x,y)}$ is the  zero measure. Then the  state of the system at  time $t\ge0$ is
represented by the random point measure
\[%\label{eq:seq-mu-def}
  \Z(t) = \sum_{j=1}^{Z(0)}\delta^+_{(B^0_j(t),D^0_j(t))}
        + \sum_{i=1}^{E(t)}\delta^+_{(B_i(t),D_i(t))}.
\]
Note that the queue length at time $t$ is given by the
total mass of the measure $\Z(t)$,
\[%\label{Z=totalMass}
	Z(t)=\langle 1,\Z(t)\rangle,
\]
where $\langle f,\mu\rangle=\int_{\overline{\BR}_+^2}fd\mu$ for a Borel
measure $\mu$ on $\overline{\BR}_+^2$ and a $\mu$-integrable function
$f:\overline{\BR}_+^2\ar\BR$.

\begin{figure}
\begin{center}
\begin{pspicture}(-.7,-4)(10,4)
		%\psframe(0,-4)(10,3)
		\psaxes[linewidth=1pt,ticks=none,arrowsize=2pt 4,
						labels=none,Ox=0,Dx=1,Oy=0,Dy=1]
		  {->}(0,-3)(0,-3)(9,3)
		\rput(5,-3.3){\small\em Abandonment}
		\rput(-.7,-.4){\parbox[c]{2cm}{\small\em
			 Service \\ completion}}
		\rput(9,-3.5){ Residual service times}
		\rput(0,3.5){ Lead times}
		\psset{radius=1.5, fillcolor=black, linecolor=black, linewidth=.8pt}
		
		  \psline[linecolor=gray,linestyle=dotted]{->}(3,3)(2.2,1.7)
		\qdisk(3,3){1.7pt}

		  \psline[linecolor=gray,linestyle=dotted]{->}(.9,.8)(.1,-.5)
		\qdisk(.9,.8){1.7pt}

		  \psline[linecolor=gray,linestyle=dotted]{->}(3,1.2)(2.2,-.1)
		  \psline[linecolor=black]{->}(3,1.2)(2.2,1.2)
		  \psline[linecolor=black]{->}(3,1.2)(3,-.1)
		  \psline[linecolor=black]{-}(0,1.2)(.2,1.2)
		  \psline[linecolor=black]{-}(3,-3)(3,-2.8)
		   \rput(-.5,1.2){$D_i(t)$}
		   \rput(3,-3.3){$B_i(t)$}
		\qdisk(3,1.2){1.7pt}

		  \psline[linecolor=gray,linestyle=dotted]{->}(6,-1.5)(5.2,-2.8)
		\qdisk(6,-1.5){1.7pt}

		  \psline[linecolor=gray,linestyle=dotted]{->}(6,0)(5.2,-1.3)
		\qdisk(6,0){1.7pt}

		  \psline[linecolor=gray,linestyle=dotted]{->}(2,-1)(1.2,-2.3)
		\qdisk(2,-1){1.7pt}

		  \psline[linecolor=gray,linestyle=dotted]{->}(5,2)(4.2,.7)
		\qdisk(5,2){1.7pt}

		\rput(3.2,.55){\small$1$}
		\rput(2.6,1.5){\small$\scriptstyle {1}/{Z(t)}$}
	\end{pspicture}
\caption{Dynamics of the measure valued process ${\mathcal Z}(\cdot)$}
\label{fig1}
\end{center}
\end{figure}
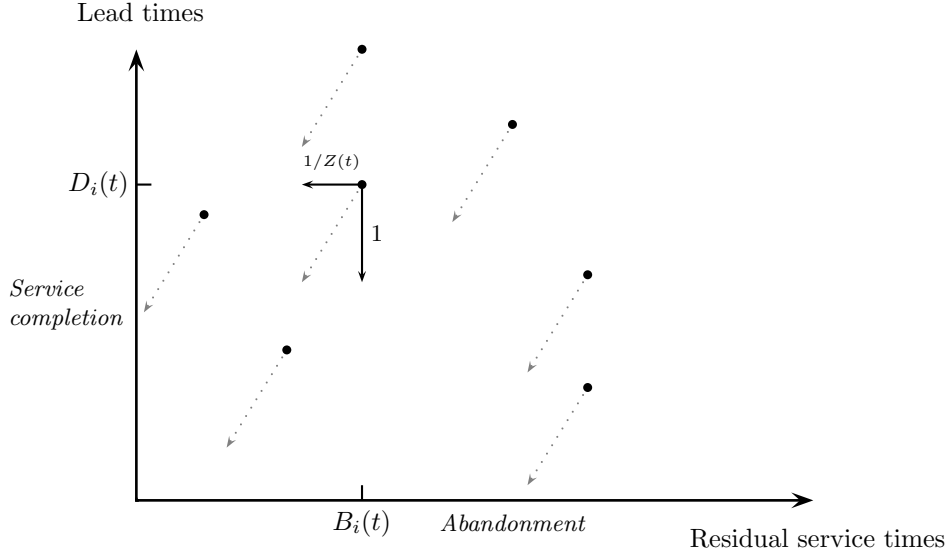

In this way, the dynamics of the  system are represented as a distribution of point masses
on $\overline{\BR}_+^2$  moving toward  the axes.  At  time $t\ge0$,  points move left  at rate
$1/Z(t)$ and down at rate $1$. (A  point with one coordinate equal to infinity will remain
that way while the  other coordinate moves.)  Point masses vanish when  hitting one of the
axes: a point mass reaching the vertical  axis corresponds to a job completing service, while a
point  mass  hitting the  horizontal  axis  represents a  job  abandoning  the queue.  See
Figure~\ref{fig1}.

Let $\M$ denote  the space of finite non-negative  Borel measures on $\overline{\BR}_+^2$,
endowed with  the topology of weak convergence:  $\zeta_n\wk\zeta$ in $\M$ if  and only if
$\langle   f,\zeta_n\rangle\ar\langle  f,\zeta\rangle$   for   all  continuous   functions
$f:\overline{\BR}_+^2\ar\BR$ (recall that $\overline{\BR}_+^2$  is compact for the induced
topology).  Let $\mathbf{D}([0,\infty),\M)$ denote the  space of {\em c\`adl\`ag} paths in
$\M$, endowed with the Skorohod $J_1$-topology.   Then, for $t\geq 0$, $\Z(t)$ is a random
element   of   $\M$   for  each   $t\ge0$,   and   $\Z\cd$   is   a  random   element   of
$\mathbf{D}([0,\infty),\M)$.

%Note that the initial measure $\Z(0)$ is a convenient way to express the initial condition.

It is clear that, given   stochastic  primitives  $(E\cd,\{B_i,D_i\}_{i=1}^\infty)$   and  the  initial
condition $\Z(0)$,  the equation \eqref{eq:seq-S} uniquely determines the processes $S\cd$,
$Z\cd$, $\Z\cd$,  and the residual service  times and lead  times. It is also  easily seen
that the  state descriptor $\Z\cd$  satisfies the following  equation: for each  Borel set
$A\in\mathscr{B}$, and all $t\ge 0$,
\begin{equation}
\label{e.dynamicEQ}
  \Z(t)(A) =
   \Z(0)\left(A+(S(t),t)\right)
   + \sum_{i=1}^{E(t)}
      1^+_A\left(B_i(t), D_i(t)\right),
\end{equation}
where $A+w=\{a+w :a\in A\}$ and $1^+_A(w) =\langle 1_A,\delta^+_{w}\rangle$. 
Note that the quantity $\Z(0)\left(A+(S(t),t)\right)$ corresponds to a shift by the
quantity $(S(t),t)$ of the initial points: indeed, if $(x,y)\in\BR_+^2$ and
$(s,t)\in\BR_+^2$, for $A\in\mathscr{B}$,    
\[
\delta_{(x,y)}(A+(s,t))=\delta_{(x-s,y-t)}(A).
\]
This equation plays a crucial r\^ole in determining fluid limits for the model. 

\subsection{A Fluid Scaling}\label{scn:scaling}
A sequence of renormalized stochastic  processes $(\Zbr(t))$ associated to the solution of
the evolution  equation~\eqref{e.dynamicEQ} is introduced. The limits  of $(\Zbr(t))$ will
give the fluid limits of this queue.

Let $\R\subset  [0,\infty)$ be a  sequence increasing to  infinity. Suppose that  for each
$r\in\R$, there is  a stochastic model as defined in  Section \ref{s.stochModel}. That is,
for  each  $r\in\R$, there are stochastic primitives  $(E^r\cd,\{B^r_i,D^r_i\}_{i=1}^\infty)$  with
associated data  $\lambda^r$ and $\vartheta^r$, and an initial  condition $\Z^r(0)$ which 
give stochastic processes $Z^r\cd$, $S^r\cd$, $\Z^r\cd$, and residual service times and lead
times $\{B^r_i\cd,D^r_i\cd\}$ and $\{B^{0r}_i\cd,D^{0r}_i\cd\}$.  Each model is defined on
a   probability   space   $(\Omega^r,\mathscr{F}^r,\BP^r)$   with   expectation   operator
$\BE^r(\cdot)$.

A fluid scaling  is applied to each model in  the sequence.  To
obtain non-trivial scaling limits, initial lead  times $\{D^r_i\}$ will be assumed to be of
order $r$.  For each $r\in\R$,  let $\breve{\vartheta}^r\in\M$ be the  probability measure
defined by
\[
\breve{\vartheta}^r(F\times G) = \vartheta^r(F\times r G),
\]
for all  Borelian subsets $F,G  \in\mathscr{B}$ with the notation $rG=\{r\cdot g:g\in
G\}$. Note that if $\{B_i^r,D_i^r\}=\{B_i,rD_i\}$ for  some sequence
$\{B_i,D_i\}$, then $\breve{\vartheta}^r$ is simply the distribution
of $(B_1,D_1)$.

For each $r\in\R$, the fluid scaled state descriptor is defined, for
$t\ge 0$, as the random measure $\Zbr(t)\in \bM$ such that
\begin{equation}
\nonumber
\Zbr(t)(F\times G) = \frac{1}{r}\Z^r(rt)(F\times rG),
\end{equation}
for all  Borelian subsets $F,G  \in\mathscr{B}$. This definition  scales lead times by  a factor
$r^{-1}$  as well.   Fluid  scaled versions  of  the remaining  processes  are defined  as
follows: for all $r\in\R$, $t\ge s\ge0$, and $i=1,\ldots,E^r(rt)$, let
\begin{xalignat*}{2}
 \Er(t) & =  \frac{1}{r}E^r (rt), & \overline{Z}^r(t) &= \frac{1}{r}Z^r(rt),\\
 \overline{S}^r(t) & = S^r(rt), & \overline{S}^r(s,t) & = \overline{S}^r(t)-\overline{S}^r(s),\\
\overline B^r_i(t)  &= B^r_i(rt), & \overline D^r_i(t) & = \frac{1}{r}D^r_i(rt).
\end{xalignat*}

The following asymptotic assumptions are needed. Let
$(\lambda,\vartheta,\zeta_0)$ be fluid model data satisfying the
assumptions of Section \ref{s.fluidModel}. Assume that as $r\ar\infty$,
\begin{align}
\label{a.u-conv-exp} 
&\Er(\cdot)\to\lambda(\cdot), \\ % \lambda^r & \lar \lambda, \\
\label{a.lawConv}         
&\breve{\vartheta}^r  \wk  \vartheta,\\
\label{a.initconv}
&\Zbr(0) \wk \zeta_0,\quad \text{in distribution},
\end{align}
in particular, Assumption~\eqref{a.u-conv-exp} implies that $\lambda^r\to\lambda$
holds.  

\subsection{Fluid model.}\label{s.fluidModel}
A deterministic fluid model  satisfying dynamic equations analogous to \eqref{e.dynamicEQ}
is introduced.  It will  be shown later that these equations can  be obtained as limits of
Equation~\eqref{e.dynamicEQ} under an appropriate scaling procedure.

\begin{comment}
For  the   sum  on   the  right   side  of
\eqref{e.dynamicEQ},  this will  require a  functional Glivenko-Cantelli  theorem  that is
uniform over  the sets  $A$ used  in the equation  (see Section  \ref{s.GC}). Thus,  it is
important  to   take  the  limit  using   a  class  $\mathcal{C}$  of   Borel  subsets  of
$\overline{\BR}_+^2$ that  is small  enough for this  uniformity to  hold, yet large  enough to
still characterize elements  of $\M$. This is accomplished in  Section \ref{s.GC} with 
\end{comment}

Let $\rho>1$, and $\zeta_0\in\M$ be a measure on $\overline{\BR}_+^2$ such that the projections
$\zeta_0(\cdot\times\overline{\BR}_+)$ and $\zeta_0(\overline{\BR}_+\times\cdot)$  are free of atoms
in $[0,\infty)$ and $z_0=\zeta_0(0,0)$ is the total mass of $\zeta_0$.

\begin{defn}\label{d.fms}
A {\em measure valued fluid  model solution} for the data $(\lambda,\vartheta,\zeta_0)$ is
a  continuous  function $\zeta\cd:[0,\infty)\ar\M$  such that
\begin{enumerate}
\item \label{prop-a} $\inf_{t>a} z(t)>0$ for all $a>0$,
\item \label{prop-b} for all $C\in\mathcal{C}$ and $t\ge0$,
\begin{equation}
\label{e.mvfms}
     \zeta(t)(C)  =\zeta_0(C+(S(0,t),t))
       +\lambda\int_0^t \vartheta(C+(S(s,t),t-s))\,{\rm d}s,
\end{equation}
where  $S(u,v)=\int_u^v{1}/{z(s)}\,{\rm d}s$ for all $v\ge u\ge0$ and $z\cd$ is the total
mass function is $z\cd=\langle 1,\zeta\cd\rangle$. The  function $z\cd$ is simply called a
{\em fluid model solution}  for $(\lambda,\vartheta,\zeta_0)$.
\end{enumerate}
\end{defn}
\noindent
Note that  $S(0,t)$ may be equal to  $+\infty$ if $\zeta_0\equiv 0$,  i.e. $z(0)=0$.  Both
right hand side terms in \eqref{e.mvfms} are  still well defined  in this case, and  the first
term equals zero.

The class of corner sets is defined as 
\[%\label{a.Cdef}
    \mathcal{C}=\left\{ [x,\infty)\times[y,\infty):x,y\in\Rp \right\}
        \cup\left\{ [x,\infty]\times[y,\infty]:x,y\in\overline{\BR}_+ \right\}.
%	\mathcal{C}=\left\{ [x,\infty]\times[y,\infty]:(x,y)\in\overline{\BR}_+^2\right\}.
\]
The sets  from the class  $\mathcal{C}$ will  be used to  describe the evolution  of fluid
model solutions. Since  each $C=[x,\infty]\times[y,\infty]\in\mathcal{C}$ is characterized
by  the  coordinates  $(x,y)$  of  its  corner,  it is  convenient  to  use  the  notation
$\mu(x,y)\stackrel{\text{def.}}{=}\mu([x,\infty]\times[y,\infty])$       for       any
$\mu\in\M$.  If $z_0>0$, for this  class of subsets Equation~\eqref{e.mvfms} can then be
rewritten as follows: for each $x,y\geq 0$ and $t\ge0$,
\begin{multline}
\label{measurefluid}
\zeta (t)(x,y)=z_0\pr{B^0>x+S(0,t), D^0>y+t}\\+\lambda
\int_0^t \pr{B>x+ S(s,t), D>y+t-s} \,{\rm d}s,
\end{multline}
Since $z(t)=\zeta(t)(0,0)$, the  fluid model solution  $z\cd$ satisfies the  following equation:
for each $t\ge0$,
\begin{equation}\label{q}
z(t)=z_0\pr{ B^0 > S(0,t); D^0 > t} + \lambda
\int_0^t \pr{ B > S(s,t); D > t-s}\,{\rm d}s.
\end{equation}

%The data $(\lambda,\vartheta,\zeta_0)$ correspond to limiting values of the data
%$\lambda$, $\vartheta$, and $\Z(0)$ used to define the stochastic model.

It will be proved  that the fluid model defined above is the  limit in distribution of the
rescaled  processes $\{\Zbr\cd:r\in\R\}$  introduced in  Section  \ref{s.stochModel}.  The
measure valued  fluid model  solution $\zeta\cd$ corresponds  to the measure  valued state
descriptor $\Z\cd$,  and the  fluid model solution  $z\cd$ is  the limit of the  queue length
process $Z\cd$.  The main result concerning the convergence of $\{\Zbr\cd:r\in\R\}$ is the
following theorem.
\begin{thm}\label{t.main}
The sequence  $\{\Zbr\cd:r\in\R\}$ is tight and each  weak limit point is  almost surely a
measure     valued     fluid     model     solution    $\zeta\cd$     for     the     data
$(\lambda,\vartheta,\zeta_0)$.  If  in addition  \eqref{a.lip}  holds, then  $\Zbr(\cdot)$
converges  in distribution,  as $r\ar\infty$,  to the  unique measure  valued  fluid model
solution $\zeta(\cdot)$.
\end{thm}
\noindent
This theorem is proved in Section~\ref{s.tight} and~\ref{s.limits}.

\subsection{Some Properties of the fluid model}\label{scn:result}
Despite the quite abstract  setting of this paper (measure valued
processes), some concrete and explicit results concerning the fluid model of the queue can
be    obtained.   Let   $(\lambda,\vartheta,\zeta_0)$    satisfy   the    assumptions   of
Section~\ref{s.fluidModel}.

The first result establishes the uniqueness of fluid model solutions under a Lipschitz
condition on the initial condition $\zeta_0$. 
\begin{thm}\label{t.lip}
Suppose there exists a finite constant $L$ such that 
\begin{equation}
\label{a.lip}
\zeta_0(A \times [y,y']) \leq L|y'-y|
\end{equation}
for all Borel sets $A$ of $\overline\BR$ and all $y'>y\geq 0$ . Then
\eqref{measurefluid} and \eqref{q} have a unique solution.  \end{thm}
\noindent

The second theorem analyzes the equilibrium of the fluid model, i.e.  the behavior at
infinity of the solution of Equation~\eqref{q}.
\begin{thm}
Suppose that $\lambda\expect{B1_{\{D=\infty\}}}<1$ and $\expect{\min
\{B,D\}}<\infty$. Then any solution $t\to z(t)$ of \eqref{q} converges
at infinity to the unique positive solution $z_\infty$ of the fixed point equation 
\[
z_\infty=\lambda     \expect{\min    \{z_\infty     B,D\}}.
\]
Moreover,  any solution  $(\zeta(t))$  of Equation~\eqref{measurefluid}  converges to  the
measure $\zeta_\infty\in\M$, defined by
\begin{multline}\label{Rully}
\zeta_\infty(x,y)=\lambda \int_0^\infty \pr{B>x+\frac{t}{z_\infty}, D>y+t}\,{\rm d}t
\\= \expect{\min \{z_\infty(B-x)^+, (D-y)^+ \}},
\end{multline}
for  $x,y\ge0$.
\end{thm}
These theorems  are proved in Section  \ref{s.analysis}.  The simple  fixed point equation
stated in  this theorem  is used  to analyze the  qualitative behavior  of the  queue, see
Section~\ref{s.applications}.  Note that  the  expression of  the  distribution of  points
$\zeta_\infty$ describing the asymptotic behavior of this queue has a simple expression in
terms of the solution $z_\infty$ of the fixed point equation.

\section{Some Properties of the Fluid Model}\label{s.analysis}
\setcounter{equation}{0} 

%The section  is organized as follows.  Section \ref{s.fluid.bound} shows that all  fluid
%model solutions are  contained in a certain compact  set, by showing that each fluid
%model solution is bounded, and by constructing  a {\em maximal} solution. The   behavior
%at  infinity   of  fluid   model  solutions   is  considered   in  Section
%\ref{s.fluid.infinity}. It turns out that the condition $\rho>1$ leads to nontrivial
%limit behavior,  capturing the balancing  effects of  overload and  impatience.  After
%that, we investigate uniqueness of fluid model solutions. Section \ref{s.fluid.unique1}
%establishes uniqueness under nonzero  conditions.  If the system starts empty  at time
%$0$, considered in  Section \ref{s.fluid.unique2},  uniqueness  is more  delicate.
%Finally, we  consider analogous properties of measure valued fluid model solutions.
%\subsection{Basic properties}\label{s.fluid.basic}

In this section some basic properties of fluid model solutions are derived.  
In what follows, let $z(\cdot)$ be an arbitrary fluid model solution, i.e. such that
\begin{equation}\tag{\ref{q}}
z(t)=z_0\pr{ B^0 > S(0,t); D^0 > t} + \lambda
\int_0^t \pr{ B > S(s,t); D > t-s}\,{\rm d}s.
\end{equation}
If $z_0>0$, define 
\[
\widetilde S(t)=\inf \{s: S(0,s)\geq t\},
\]
since $z(t)\leq \lambda t+z_0$, $S(0,t)\rightarrow\infty$ as
$t\rightarrow\infty$, implying that $\widetilde S(t)$ is well defined for
all $t$.  In addition, $\widetilde S(t)<\infty$ for all $t$ if $z_0>0$,
which follows from property \eqref{prop-a} and continuity of the fluid
model solution $z(\cdot)$. 

Define  $\widetilde z(t)=z(\widetilde S(t))$, then $(\widetilde z(t))$  satisfies the  
equation 
\begin{multline}\label{tildeq}
\widetilde z(t) =  z_0 \pr{ B^0 \geq t;D^0 \geq \widetilde S(t)}\\+\lambda
\int_0^t \widetilde z(u) \pr{B\geq t- u; D \geq \widetilde S(t)-\widetilde S(u)}\,{\rm d}u.
\end{multline}

\iffalse
If $t^*=\inf \{t: z(t)=0 <\infty\}$, then $S(0,t)=\infty$ for $t\geq
t^*$; in this case $\widetilde S(t)\rightarrow t^*$ as $t\rightarrow\infty$.
We  address the question when $t^*<\infty$ and show that $z(t)=0$ for
$t>t^*$ later on. For the analysis of the time-changed process $\widetilde
z(t)$ we do not need to distinct between the cases $t^*<\infty$ and
$t^*=\infty$.  An important distinction we do need to make is to whether
$z_0=0$ or not. If $z_0=0$, the time change $t\rightarrow \widetilde S(t)$
is not possible, since $S(0,t)$ can be $\infty$ in this case.  
\fi

We next introduce the concept of a {\em shifted} fluid model solution:
For $t_0>0$ define $z(t_0,t)=z(t_0+t)$, and define
$S(t_0,u,v)= \int_u^v 1/z(t_0,r) {\rm d}r$.

\begin{property}
If $z(\cdot)$ is a solution of (\ref{q}), then $z(t_0,\cdot)=z(t_0+
\cdot)$ satisfies
\[
z(t_0,t)= \zeta(t_0)(S(t_0,0,t),t)+ \lambda \int_0^t \pr{B\geq S(t_0,s,t); D\geq t-s}
\,{\rm d}s, \vspace{1cm} t\geq 0.
\]
\end{property}
\begin{proof}
Note first that by definition, and since $S(t_0,0,t)=S(t_0,t_0+t)$,
\begin{align*}
\zeta& (t_0)(S(t_0,0,t),t) = z_0\pr{B^0\geq S(t_0,0,t)+S(0,t_0),D^0\geq t_0+t} 
\\& \qquad \qquad+
\lambda
\int_0^{t_0} \pr{B\geq  S(t_0-s,t_0; D \geq  s+t)+S(t_0,0,t)} {\rm d}s 
\\&= z_0 \pr{B^0\geq  S(0,t_0+t) ;D^0 \geq  t_0+t}\\
&  \qquad \qquad+
\lambda \int_0^{t_0} \pr{ B\geq  S(t_0-s,t_0+t); D \geq  s+t} {\rm d}s.\\
&= z_0 \pr{B^0\geq  S(0,t_0+t) ;D^0 \geq  t_0+t}\\
& \qquad \qquad +
\lambda \int_0^{t_0} \pr{ B\geq  S(s,t_0+t); D \geq  t_0+t-s} {\rm d}s.
\end{align*}
We apply this expression as follows:
\begin{align*}
z(t_0,t) &= z(t_0+t) \\
&= z_0 \pr{ B^0\geq  S(0,t_0+t); D^0\geq t_0+t}\\
&\qquad\qquad\qquad +\lambda \int_0^{t_0+t} \pr{ B\geq  S(s, t_0+t); D\geq t_0+t-s} {\rm d}s\\
&= \zeta (t_0)(S(t_0,0,t),t)+
 \lambda \int_{t_0}^{t_0+t} \pr{B\geq  S(s, t_0+t); D\geq t_0+t-s} {\rm d}s.
\end{align*}
The change of variables $y=s-t_0$ and the identity
$S(t_0+y,t_0+t)=S(t_0,y,t)$ give the result.
\end{proof}

The next proposition shows  that continuity of fluid
model solutions is a consequence of properties \eqref{prop-a} and
\eqref{prop-b}.

\begin{lem}
\label{continuity}
Let the distribution of $(B^0,D^0)$ be free of atoms. Then
any solution $z(t),t\geq 0$ to (\ref{q}) satisfying $\inf_{t>a} z(t)>0$
for all $a>0$ is continuous.
\end{lem}

\begin{proof}
The function $t\to S(0,t)$ is continuous and so is $t\to \pr{ B^0 > S(0,t); D^0 > t}$
since the distribution of $(B^0,D^0)$ has no atom. The first term of the right hand side
of Equation~\eqref{q} is a continuous function of $t$. Concerning the second term of
Equation~\eqref{q}, by monotonicity the integrand $t\to \pr{ B > S(s,t); D > t-s} $ is
continuous almost everywhere on $\BR_+$, hence its integral is a continuous function of
$t$ by Lebesgue's Theorem. The lemma is proved. 
\end{proof}

\begin{comment}
Let $\varepsilon >0$. It is straightforward to show that
$|z(t+\varepsilon)-z(t)|\leq I+II$, with
$$I=z_0\pr{S(0,t)\leq B^0< S(0,t+\varepsilon);t\leq D^0< t+\varepsilon},$$
and
$$
II \leq \varepsilon + \int_0^t \pr{S(s,t)\leq B<S(s,t+\varepsilon)} {\rm d}s.
$$
Since $(B^0,D^0)$ is free of atoms, and $S(0,t)$ is continuous in $t$ and
finite if $z_0>0$, we have that $I\rightarrow 0$ as $\varepsilon\downarrow
0$.  The second term converges to $ \int_0^t \pr{ B=S(s,t)} {\rm d}s
$ which is $0$ since for every $s>0$, $S(s,t)$ is finite, strictly
increasing and continuous in $s$, and the distribution function of $B$
only has countably many discontinuities.  
\end{comment}

\subsection{A Maximal Solution}
\label{s.fluid.bound}
An important monotonicity property of fluid model solutions is proved in this section. 
\begin{prop}
\label{upperbound}
If $\lambda \expect{B 1_{\{D=\infty\}}}<1$ and $\expect{\min\{B,D\}}<\infty$,
then any fluid model solution is bounded.
\end{prop}

\begin{proof}
Note first that, since $\expect{\min \{B,D\}}<\infty$, also $\expect{\min
\{B,aD\}}<\infty$ for every $a\in [0,\infty)$. Define $\|z\|_t=\sup_{0\leq
u\leq t} z(u)$. Note that
$\|z\|_t\leq z_0+\lambda t<\infty$. Fix $t$ and let $u\in [0,t]$. Since
$S(u,s) \geq (u-s)/\|z\|_t$, it holds that
\[
z(u)\leq z_0+\lambda \int_0^u \pr{D\geq u-s; \|z\|_t B \geq  u-s} {\rm d}s \leq z_0+\lambda \expect{\min\{D, \|z\|_t B\}}
\]
which is finite since $\expect{\min \{D,B\}}<\infty$.
By taking the supremum over $u\in[0,t]$ and by dividing both sides by $\|z\|_t$ one
obtains the relation
\[
1 \leq z_0/\|z\|_t + \lambda\expect{\min\{D/\|z\|_t, B\}},
\]
If $\|z\|_t\rightarrow\infty$ then, by monotone convergence, one gets the inequality
\[
1 \leq \lambda\expect{B 1_{\{D=\infty\}}},
\]
which contradicts the assumption $ \lambda \expect{B 1_{\{D=\infty\}}}<1$.
We conclude that $\|z\|_t$ converges to some finite constant $M$ which implies
the assertion.
\end{proof}

\begin{prop}[Maximal Fluid Solution]\label{maxsol}
For $z_0>0$, there exists a fluid model solution $z^*(\cdot)$ starting form $z_0$ which is
maximal, i.e. for any fluid model solution $z(\cdot)$ such that $z(0)=z_0$, the relation $z(t)\leq z^*(t)$ holds for all $t\geq 0$.
\end{prop}
\begin{proof}
To define $z^*(\cdot)$ we first define a sequence of functions
$z^n(\cdot), n\geq 0$, by $z^0(t)=z_0+\lambda t$,
$S^n(u,v)=\int_u^v (1/z^n(r)) {\rm d}r$ and
\[
z^{n+1}(t)= z_0 \pr{B^0 \geq  S^n(0, t); D^0\geq  t}+\lambda
\int_0^t \pr{B\geq  S^n(t-s,t); D\geq s} {\rm d}s.
\]
We show that $z^{n+1}(t)\leq z^n(t)$ by induction.
The inequality $z^1(t)\leq z^0(t)$ is trivial.
Suppose now that $z^n(t)\leq z^{n-1}(t)$. Then $S^n(u,v)\geq
S^{n-1}(u,v)$, and, using the fact that tail probabilities are non-increasing,
\begin{align*}
z^{n+1}(t) &= z_0 \pr{B^0 \geq  S^n(0, t); D^0\geq  t}+\lambda
\int_0^t \pr{D\geq s; B\geq  S^n(t-s,t) } {\rm d}s\\
&\leq  z_0 \pr{B^0 \geq  S^{n-1}(0, t); D^0\geq  t}+\lambda
\int_0^t \pr{B\geq  S^{n-1}(t-s,t); D\geq s} {\rm d}s,
\end{align*}
which equals $z^{n}(t)$.

Since $z^n(t)$ is decreasing in $n$ and non-negative for all $n$ there
exists a function $z^*(t)$ such that
$z^*(t)=\lim_{n\rightarrow\infty} z^n(t)$.
By the definition of $z^n(t)$, we see that $z^*(t)$ satisfies (\ref{q}).

Furthermore, we have $z(t)\leq z^*(t)$ for any given fluid model
solution $z(\cdot)$. This is true because $z(t)\leq z^0(t)$, and
using an inductive argument as above, $z(t)\leq z^n(t)$ for every $n$.
Since we know that at least one fluid model solution exists, it follows
that $\inf_{t>a} z^*(t)>0$ for every $a>0$. By Lemma~\ref{continuity},
it follows that $z^*(t)>0$ is continuous. We conclude that $z^*(\cdot)$
is indeed a fluid model solution.
\end{proof}

\subsection{Convergence of fluid model solutions}
\label{s.fluid.infinity}

In this subsection we show the convergence of fluid model solutions
to a non-trivial constant $z_\infty$ as $t\rightarrow\infty$.
%As a preliminary result, we first give the following proposition.

\begin{prop}
If $\lambda\expect{B1_{\{D=\infty\}}}{<}1$, $\expect{\min\{B,D\}}{<}\infty$ and
$\rho{>}1$, the equation 
\begin{equation}
\label{fixedpointeq}
z_\infty=\lambda \expect{\min\{z_\infty B,D\}},
\end{equation}
has a unique solution in $(0,\infty)$.
\end{prop}

\begin{proof} 
The  function $f:  a\to \lambda  \expect{\min\{B,aD\}}$ is  non-decreasing and  concave on
$[0,+\infty)$,  note that  $f(a)= \lambda  \expect{\min\{B,aD\}1_{\{D<+\infty\}}}+ \lambda
\expect{B1_{\{D=+\infty\}}}$        for       $a>0$,        therefore       $f(0+)=\lambda
\expect{B1_{\{D=+\infty\}}}{<}1$  and $f(a)$ converges  to $\lambda\expect{B}{>}1$  as $a$
goes to infinity.  By continuity of $f$, there exists $a_0$, $0{<}a_0{<}+\infty$ such that
$f(a)=1$. The concavity and the monotonicity  imply that such a $a_0$ is unique, otherwise
$f$  should be  constant  equal to  $1$  after $a_0$,  but this  is  impossible since  $f$
converges to $\lambda\expect{B}{>}1$ at infinity.  The quantity $1/a_0$ is then the unique
solution of Equation~\eqref{fixedpointeq}.
\end{proof}

\noindent
We are now ready to present the  main result of this subsection, concerning the asymptotic
behavior of any fluid model solution $(z(t))$ as $t$ goes to infinity.

\begin{thm}
If $z(\cdot)$ is a fluid model solution, under the conditions
\[
\lambda\expect{B1_{\{D=\infty\}}}<1, \expect{\min\{B,D\}}<\infty \text{ and }\rho>1,
\]
then $z(t){\to} z_\infty$ as $t{\to}+\infty$,
where $z_\infty$ the unique positive solution of the fixed point equation (\ref{fixedpointeq}).
\end{thm}

\begin{proof}
It suffices to show  $\bar z=\limsup_{t\rightarrow\infty} z^*(t)\leq z_\infty$ and 
$\underline{z}=\liminf_{t\rightarrow\infty} z^*(t)\geq z_\infty.$
We start with the former. We know that $\bar z<\infty$ from Proposition \ref{upperbound}.
For any $\varepsilon>0$ there exists a $t_\varepsilon$ such
that $z(t) \leq \bar z +\varepsilon$. We see that, for $t>t_\varepsilon$,
\begin{multline*}
z(t) \leq z_0 \pr{B^0\geq S(0,t); D^0\geq t} \\+
\lambda \int_0^{t_\varepsilon} \pr{B\geq S(s,t);D\geq  t-s}\,{\rm d}s+\lambda \int_0^{t-t_\varepsilon}
\pr{\min \{ D, (\bar z+\varepsilon) B\} > s }\,{\rm d}s.
\end{multline*}
Taking the $\limsup$ on both sides, and noting that $S(s,t)\rightarrow\infty$ for any
$s\geq 0$, we obtain that
\[
\bar z \leq \expect{ \min \{ D, (\bar z+\varepsilon) B\} }.
\]
The result is valid for every $\varepsilon>0$. By letting $\varepsilon\downarrow 0$,
we obtain that $\bar z\leq z_\infty$.
The lower bound follows by a similar argument after first noting that
$\underline{z}>0$ since $z(\cdot)$ is a fluid model solution.
\end{proof}

\iffalse
We conclude by investigating the large time behavior of
$z(t)$ in the case $\rho\leq 1$.

\begin{prop}
Let $z(t)$ be a solution of $(\ref{q})$.
If either $\rho<1$ or $\rho=1$ and $\pr{D>aB}<1$ for all $a>0$,
then $z(t)\rightarrow 0$.
\end{prop}

\begin{proof}
As in the previous result we have $\limsup_{t\rightarrow\infty} z(t)\leq z^*$,
with  $z^*=\lim_{n\rightarrow\infty} z^n$ defined as above. $z^*$ is the largest possible solution of
the equation $\lambda \expect{\min \{D,Bz^*\}}=z^*$. If $z^*>0$ this would imply that both
$\rho =1$ and $\pr{D>z^*B}=1$, contradicting our assumptions. \end{proof}

Note that for $\rho=1$, the additional condition on the joint distribution of $B$ and $D$ is necessary.
For example, if $D=\infty$ it holds that $t^*=\infty$ whenever $z(0)>0$, see \cite{GrPuWi2002}.
\fi

\subsection{Uniqueness of fluid model solutions under non-zero conditions}
\label{s.fluid.unique1}
The uniqueness of fluid model solutions is, in general, difficult to determine.  If one
looks at the time-changed version (\ref{tildeq}), and take $\lambda=0$, one gets an
ODE. Uniqueness of solutions to such ODE's can usually only be established by reducing it
to some special case or to assume some kind of Lipschitz condition. If $D=\infty$, then
(\ref{tildeq}) reduces to a renewal equation for which uniqueness is known to
hold. Unfortunately, this reduction is not possible in general, which lead us to use a
Lipschitz condition on the distribution function of $(B^0,D^0)$. It is not
necessary to assume regularity conditions on the distribution of $(B,D)$.  We shall give a
direct proof of uniqueness; for related results in the functional analysis literature, we
refer to Chapter 2 of Hale \& Verduyn Lunel \cite{HV93}.

\begin{thm}
\label{lip}
Suppose $z_0>0$ and $F_0(x,y)=\pr{B^0\geq x;D^0\geq y}$
is Lipschitz continuous in $Y$, i.e. there is a constant $L$ such that
for any $x,y,y'$,
\[
|F_0(x,y)-F_0(x,y')|\leq L|y-y'|.
\]
Then (\ref{q}) has a unique solution.
\end{thm}

%To prove Theorem (\ref{lip}), we introduce some additional notation.
\noindent
Defining $\widetilde \zeta(t)(u,v)=\zeta (\widetilde S(t))(u,v)$, it can easily be shown that
\begin{multline}
\label{quv}
\widetilde \zeta(t)(u,v) = z_0 \pr{B^0\geq u+t;D^0\geq v+\widetilde S(t)}\\+\lambda \int_0^t \widetilde z(s)
\pr{B\geq  u+(t-s); D\geq  v + \widetilde S(t)-\widetilde S(s)}\,{\rm d}s.
\end{multline}
It  is  clear  that for  any  $u,v$,  $\widetilde  \zeta(t)(u,v),t\geq 0$,  is  completely
determined  by $\widetilde z(t),t\geq  0$ and  the initial  measure.  Thus,  uniqueness of
$z(t)$ on an interval $A$ carries over to uniqueness of $\zeta(t)(u,v)$ on $A$.

The idea of  the proof is simple: we  take a suitable constant $a>0$ and  prove first that
uniqueness  holds for  $\widetilde  z(t)$  for $[0,a]$.   As  discussed above,  uniqueness
carries over  to $\widetilde \zeta(t)(u,v)$ for  $t\in [0,a]$. Using this  and the shifted
fluid model equation  given by Property 5.1, we prove uniqueness  for $\widetilde z(t)$ on
the  interval $[a,2a]$,  and so  forth.   This iterative  procedure works  if the  measure
$\widetilde  \zeta(t)(u,v)$ is  Lipschitz for  any $0<t<T$.   This is  the content  of the
following lemma.

\begin{lem}
\label{ll}
For any $x,y,y'$ and for any $t$ we have
\[
|\widetilde \zeta(t)(x,y) - \widetilde \zeta(t)(x,y')| \leq (z_0L+\lambda)|y-y'|.
\]
\end{lem}

\begin{proof}
We may take $y,y'$ such that $y\leq y'$.
 From (\ref{quv}) we obtain
\begin{multline*}
|\widetilde \zeta(t)(x,y) - \widetilde \zeta(t)(x,y')| \leq z_0L|y'-y|\\+
\lambda \int_0^t z(s) \pr{\widetilde S(t)-\widetilde S(s) + y \leq D \leq \widetilde S(t)-\widetilde S(s)+y'} {\rm d}s.
\end{multline*}
Noting that $\widetilde z(s){\rm d}s={\rm d}\widetilde S(s)$ we can rewrite this into
\beeq
\label{broembrom}
z_0L|y'-y|+\lambda \int_0^{\widetilde S(t)} \pr{r + y \leq D \leq r + y'} {\rm d}s
\eneq
Noting that, for any $\delta>0$
\[ %\label{deltabound}
\lambda \int_0^\infty \pr{y\leq D< y+\delta} {\rm d}y\leq \delta,
\]
we see that (\ref{broembrom}) can be upper bounded by $(z_0L+\lambda)|y-y'|$.
\end{proof}

\noindent
\begin{proof}[Proof of  Theorem \ref{lip}.]  By the one-to-one  correspondence between solutions of
(\ref{q})  and  (\ref{tildeq}), it  suffices  to show  that  (\ref{tildeq})  has a  unique
solution.   Define $a=  1/(2(z_0L+4\lambda))$.  We  first show  that (\ref{tildeq})  has a
unique solution on the interval $[0,a]$.  For that, suppose that there exist two different
solutions   $\widetilde  z(t),  0\leq   t\leq  a$   and  $h(t),   0\leq  t\leq   a$.   Set
$\varepsilon=\sup_{0\leq t\leq a}|\widetilde z(t)-h(t)|$.

Note that for any $0\leq s<t\leq a$,
$$|H(t)-H(s)-(\widetilde S(t)-\widetilde S(s))|\leq \varepsilon a.$$

Using (\ref{tildeq}) for both $z$ and $h$, together with the Lipschitz
assumption, we obtain, after some simple estimates,
\begin{align*}
&|\widetilde z(t)-h(t)| \leq z_0\left|\pr{B^0\geq t;D^0\geq  \widetilde S(t)}-\pr{B^0\geq  t;D^0\geq H(t)} \right|\\
&{+}\lambda \int_0^t |\widetilde z(s)\pr{B\geq t{-}s; D\geq \widetilde S(t){-}\widetilde S(s)}{-}h(s)\pr{B\geq t{-}s; D\geq H(t){-}H(s)}|{\rm d}s.
\end{align*}
The first term is bounded by $z_0La\varepsilon$. The second term is bounded by
\begin{align*}
&\lambda \int_0^t|\widetilde z(s)-h(s)|\, {\rm d}s \;+\\
&\lambda \int_0^t z(s)\left|\pr{B\geq t-s; D\geq \widetilde S(t)-\widetilde S(s)}{-}\pr{B\geq t-s; D\geq H(t)-H(s)}\right|{\rm d}s.
\end{align*}
Call these terms $IIa$ and $IIb$. We have $IIa\leq \lambda \varepsilon a$.
To bound $IIb$, we use the bound
\begin{multline*}
\left|\pr{B\geq t-s; D\geq \widetilde S(t)-\widetilde S(s)}-\pr{B\geq t-s; D\geq H(t)-H(s)}\right|  \\
\leq \pr{\widetilde S(t)-\widetilde S(s)-\varepsilon a < D\leq \widetilde S(t)-\widetilde S(s)+\varepsilon a}
\end{multline*}
to obtain (after a change of variable $r=\widetilde S(s)$)
\[
IIb \leq \lambda \int_0^{S(t)} \pr{r-\varepsilon a <D\leq r+\varepsilon a}\,{\rm d}s \leq \lambda 2\varepsilon a.
\]
Putting everything together, we see that for $t\in [0,a]$,
\[
|\widetilde z(t)-h(t)| \leq z_0La\varepsilon +  3\lambda \varepsilon a \leq \varepsilon/2,
\]
which implies that $\varepsilon=0$, i.e. that $\widetilde z(t)$ and $h(t)$ coincide on $[0,a]$.
Hence Equation~\eqref{tildeq} has a unique solution in the interval $[0,a]$. \\

Suppose now that (\ref{tildeq}) has a unique solution on $[0,ka]$ for some $k\geq 1$,
and consider the equation
\begin{multline}\label{qka}
\widetilde z(ka,t)= \zeta (ka)(t,\widetilde S(ka,t))\\
+\lambda \int_0^t \widetilde z(ka,s) \pr{B \geq  t-s ; D \geq  \widetilde S (ka,t)-
\widetilde S(ka,s)} {\rm d}s.
\end{multline}
We now show that this equation has a unique solution on $[0,a]$, implying that there exist
a unique  solution of  (\ref{tildeq}) on the  interval $[0,(k+1)a]$.   Suppose $\widetilde
z(ka,t)$  and  $h(t)$  both  satisfy (\ref{qka}),  and  set  $\varepsilon=\sup_{t\in[0,a]}
|z(ka,t)-h(t)|$.  As before, we have
$$|(\widetilde S(ka,t)-\widetilde S(ka,s)-(H(t)-H(s))|< (t-s)\varepsilon< a\varepsilon.$$
Using this, we get as before (using now the Lemma for the first term)
that
\[
|\widetilde z(ka,t)-h(t)| \leq (z_0L+\lambda )a\varepsilon + 3\lambda a\varepsilon \leq \varepsilon/2,
\]
which implies that $\varepsilon=0$, and that uniqueness of solutions of (\ref{tildeq}) holds
on the interval $[0,(k+1)a]$. Iterating this argument completes uniqueness
for all $t$.
\end{proof}

\subsection{Uniqueness starting from zero}
\label{s.fluid.unique2}

%In  this subsection  we  suppose that  uniqueness  has been  settled  for nonzero  initial
%conditions, and  show how to use this  to get uniqueness under  zero conditions. Lipschitz
% assumptions are not necessary in this part  of the proof. 
The result in this subsection can be seen as an  extension  of  a  result  of
\cite{PuStoWi2004},  who  considered  the  case PS queue without impatience, i.e. with
$D\equiv +\infty$.

\begin{thm}
Let $\varepsilon>0$.
Suppose that $(B,D)$ and a non-increasing function $F_\varepsilon(x,y)$, with
$0\leq F_\varepsilon(x,y)\leq \lambda \varepsilon$ are such that
\[
z_\varepsilon(t)= F_\varepsilon(S_\varepsilon(0,t),t)+\lambda\int_0^t \pr{B\geq  S_\varepsilon(t-s,t);D\geq  t-s}
{\rm d}s
\]
has a unique solution $z_\varepsilon(t)$ satisfying $\inf_{t>a} z_\varepsilon (t)>0$ for $a>0$. Then
$z_\varepsilon(t)\rightarrow z_0^*(t)$ as $\varepsilon\downarrow 0$.
\end{thm}
\begin{proof}
As in the construction of the maximal fluid solution, 
$z_\varepsilon (\cdot)$ can be defined as the pointwise limit $\lim_{n\rightarrow\infty} z^n_\varepsilon (\cdot)$ with $z^n_\varepsilon (\cdot)$ recursively defined by
$z_\varepsilon^0=\varepsilon +\lambda t$ and
$$z_\varepsilon^{n+1}(t)=F_\varepsilon(S_\varepsilon^n(0,t),t)+
\lambda\int_0^t \pr{D\geq s; B\geq  S_\varepsilon^n(t-s,t)} {\rm d}s.
$$
 From this construction it can be easily shown that $z_\varepsilon^n(t)$
is decreasing in $n$, and that $z_\varepsilon^n(t)\geq z_0^n(t)$.
Since also $z^*_\varepsilon(t)\leq z_\varepsilon^n(t)$,
we see that
\[
\limsup_{\varepsilon\downarrow 0} z_\varepsilon(t)\leq
\limsup_{\varepsilon\downarrow 0} z_\varepsilon^n(t)= z_0^n(t).
\]
Since this holds for any $n$, and $z_0^n(t)\rightarrow z_0^*(t)$,
we can let $n\rightarrow\infty$ to obtain
\[
\limsup_{\varepsilon\downarrow 0} z_\varepsilon(t)\leq z_0^*(t).
\]
To prove the other bound, we observe by induction
and the properties of $F_\varepsilon$ that
$z_0^n(t)\leq z_\varepsilon^n(t)$
for every $n\geq 0$.
Consequently,
$$
z_0^*(t)=\lim_{n\rightarrow\infty} z_0^n(t)\leq
\limsup_{n\rightarrow\infty} z_\varepsilon^n(t)=z_\varepsilon(t).
$$
We conclude that $z_\varepsilon(t)\geq z_0(t)$ for every $\varepsilon>0$,
which implies the lower limit and the convergence
$z_\varepsilon(t)\rightarrow z(t)$.
\end{proof}

Uniqueness of fluid model solutions starting from $0$ is now a simple
corollary.

\begin{cor}
Suppose that
$z_0=0$.
Then (\ref{tildeq}) has a unique solution.
\end{cor}

\begin{proof}
Let $z(\cdot)$ be a fluid model solution.
Define $z_\varepsilon(t)=z(t+\varepsilon)$.
Given $z(s), 0\leq s\leq \varepsilon$, $z_\varepsilon(\cdot)$ satisfies the equation
\[
z_\varepsilon(t)=F_\varepsilon(S_\varepsilon(0,t),t)+
\lambda\int_0^t \pr{D\geq  s; B\geq  S_\varepsilon(t-s,t)} {\rm d}s.
\]
Here (with obvious notation)
\[
F_\varepsilon (x,y)=\int_0^\varepsilon \pr{D\geq s+y; B\geq x+
\int_{\varepsilon-s}^\varepsilon \frac 1{z(u)}{\rm d}u}
{\rm d}s.
\]
We see that $F_\varepsilon$ is globally Lipschitz in the second coordinate (with Lipschitz
constant 1).  Consequently, the above equation  has a unique fluid model solution in terms
of $F_\varepsilon$ so that $z_\varepsilon  (\cdot)$ is uniquely determined by $z(t), 0\leq
t\leq \varepsilon$.   Since $F_\varepsilon(x,y)\leq\lambda  \varepsilon$, we see  from the
previous    theorem    that     $z_\varepsilon(t)\rightarrow    z_0^*(t)$.     But    also
$z_\varepsilon(t)=z(t+\varepsilon)\rightarrow  z(t)$,  since  $z(t)$   is  continuous.   We  conclude  that
$z(t)=z_0^*(t)$, which implies uniqueness.
\end{proof}

\subsection{Analysis of the measure-valued fluid model}\label{s.fluid.measure}

%This  subsection  states  some   properties  for  the  measure-valued  process
%$\zeta\cd$ associated to  a fluid  model solution  $z\cd$.  

For any  Borel set $F$ of $\overline{\BR}^2_+$, the measure valued function $\zeta\cd$ 
satisfies the equation
\begin{equation}
\label{gen.fluid.eq}
\zeta(t)(F)=\zeta_0(F+(S(0,t),t))+\lambda \int_0^t \vartheta(F+ (S(s,t),t-s)) {\rm d}s.
\end{equation}
The properties, which
are analogues of properties of $z\cd$,  are gathered in the following
theorem:

\begin{thm}
Let $\zeta\cd$ be a solution of (\ref{gen.fluid.eq}).
\begin{enumerate}
\item Suppose $\rho>1, \expect{\min\{B,D\}}<\infty$ and $\lambda {\mathbf E}[B 1_{\{D=\infty\}}]<1$.
As $t\rightarrow\infty$, $\zeta(t)$ converges to the limiting measure $\zeta_\infty$
defined by
\[
\zeta_\infty ([x,\infty]\times [y,\infty]) = \lambda \int_0^\infty \pr{B-x\geq s/z_\infty,D-y\geq s}{\rm d}s,
\]
where $z_\infty$ is the unique solution of the fixed point equation~\eqref{fixedpointeq}.
\item If Condition~(\ref{a.lip}) of Theorem~\ref{t.lip} holds, then
  Equation~\eqref{gen.fluid.eq} has a unique solution. 
\end{enumerate}
\end{thm}
\begin{proof}
We know that $z(t)\rightarrow z_\infty$ as $t\rightarrow\infty$. Consequently,
$S(t-t_0,t)\rightarrow t_0/z_\infty$ for 
every $t_0>0$. Write for any Borel set $F$,
\begin{align*}
\zeta(t)(F) &= \zeta(0)(F+(S(0,t),t)) + \lambda\int_0^{t_0} \vartheta (F+(S(t-s,t),s)) {\rm d}s\\
& \hspace{0.5cm} +\lambda\int_{t_0}^t \vartheta (F+(S(t-s,t),s)) {\rm d}s.
\end{align*}
Number the three terms on the right hand side as $I,II,III$. By shifting time if necessary, we may assume that $z(0)>0$.
The first term converges to $0$. Since $z(0)>0$, there exists an $\eta>0$ such that $z(t)\geq\eta$
for all $t\geq 0$. This implies that $S(t-s,t)\leq s/\eta$.
Consequently, since $F\subseteq \overline \BR_2^+$,
\[
III\leq \lambda \int_{t_0}^t {\bf P}(B\geq s/\eta; D\geq  s) {\rm d}s.
\]
From this bound, it follows that $III\rightarrow 0$ as $t_0\rightarrow\infty$.
Since $\vartheta$ only has countably many discontinuities, and
$S(t-s,t)\rightarrow s/z$ on $[0,t_0]$ we have that
\[
II\rightarrow \lambda\int_0^{t_0}\vartheta (F+(s/z,s)) {\rm d}s.
\]
Taking $t\rightarrow\infty$ and then $t_0\rightarrow\infty$ yield
the first statement of the theorem.

To prove the second statement,
note that (\ref{gen.fluid.eq}) has a unique solution for $F=\overline \BR_2^+$,
which uniquely determines $S(s,t)$ for all $s,t$ with $s\leq t$. Since $\zeta$ is completely determined
by $\zeta_0$, and $S(s,t)$, uniqueness follows.
\end{proof}

\section{Applications}\label{s.applications}
\setcounter{equation}{0}
In this section we analyze a number of quantitative properties of the
fluid model equation (\ref{q}). In particular, we investigate
the fixed point equation
\beeq
\label{fixedpointeq2}
z_\infty=\lambda \expect{\min \{z_\infty B,D\}}.
\eneq
We treat a number of examples which allow for explicit
computations, and also obtain a number of stochastic ordering
results.
In addition, we investigate the time-dependent
behavior of $z(t)$ for exponentially distributed lead times.

We first give a
 heuristic interpretation of  Equation (\ref{fixedpointeq2}):
Let $Z^r$ denote the steady-state number of customers in the system.
Furthermore, let $V^r(B)$ be the sojourn time of a customer
if the customer never reneges. Then the actual sojourn time is
given by $\min \{V^r(B), D^r \}$, and from Little's law we get
\beeq
\label{little}
\expect{Z^r} = \lambda \expect{\min \{V^r(B), D^r\}}.
\eneq
Divide both sides of $(\ref{little})$ by $r$.
Since we observe the system in steady state at time $0$, the number of
customers hardly changes and by the snapshot principle we conclude
that $V^r= Z^r B + {\rm o}(r)$. Furthermore, we have $D^r=Dr$.
Noting that $Z^r/r\rightarrow z$ then gives (\ref{fixedpointeq2})
after dividing both sides of (\ref{little}) by $r$ and letting
$r\rightarrow\infty$.

Apart from the mean queue length $z$, we are also interested in
the long term fraction of customers that leave the system successfully.
Denote this fraction by $P_s$. It is clear that
\[%\label{succes}
P_s=\pr{D>z_\infty B}.
\]

The following remarkable property, which simply follows from
the fixed-point equation (\ref{fixedpointeq2}),
shows that the performance of the system
does not depend on the average of $D$.

\begin{property} %\label{invarianceproperty}
Consider two systems numbered by $1$ and $2$ such that $(B_2,D_2) \equiv (B_1, aD_1)$ for some $a>0$,
and such that $\lambda_1=\lambda_2$.
Then (with obvious notation) we have
\[
z_{2,\infty} = a z_{1,\infty}, \quad P_{s,2} = P_{s,1}.
\]
\end{property}

We now proceed by analyzing a number of special cases.
In Section 4.1, we assume a strong form of dependence. Section 4.2
assumes that $B$ and $D$ are independent. We give a remarkably simple expression
for $z(t)$ in the case that $D$ has an exponential distribution.
Finally, Section 4.3 considers an example which can be used as a flow level
model for the integration of elastic and streaming traffic.

\subsection{Completely dependent lead times}

Consider first the case where $D=\Theta B$, with $\Theta>0$ (independent
of $B$) reflecting the average service rate expected by a customer.
In this case, the performance measures can be determined from the equations
(recall that $\rho=\alpha \expect{B}>1$)
\[
z_\infty = \rho \expect{\min \{\Theta, z_\infty \}}, \quad
P_s = \pr{\Theta>z_\infty }.
\]

Some specific examples:

\begin{itemize}
\item {\em $\Theta$ single-valued. } If we assume that $\Theta=\theta$,
then $z_\infty =\rho \min \{\theta, z_\infty \}$, which implies that $z_\infty =\rho \theta$
since $\rho>1$. From this, it follows that all customers leave the system
impatiently: $P_s= \pr{\theta > \rho\theta}=0$.
Observe that when a customer
leaves the system,  a fraction $1/\rho$ of his service time
has been processed.
\item {\em $\Theta$ two-valued. } From the previous example, it is clear
that the system can only get some work done if some customers are more patient
than others. In this example we assume that $\Theta$ equals $\theta_1$
with probability $p$ and $\theta_2$ with probability $1-p$.
Take $\theta_2>\theta_1$. Equation (\ref{fixedpointeq2}) now simplifies to
\[%\label{twovalued}
z_\infty =\rho p \min \{z_\infty , \theta_1\} + \rho(1-p) \min \{z_\infty ,\theta_2\}.
\]
 From this equation and the properties $\theta_2>\theta_1, \rho>1$
it follows that $z_\infty > \theta_1$.
Furthermore, $z_\infty > \theta_2$ holds if and only if the equation
\[
z_\infty =\rho p \theta_1 + \rho(1-p) z_\infty
\]
has a non-negative solution, which is the case if and only if
$\rho(1-p)<1$ (i.e. when the most patient customers cannot saturate
the system alone).
In this case we have
\[
z_\infty  = \frac{\rho p \theta_1}{1-\rho(1-p)}<\theta_2.
\]
If the last inequality is not valid or if $\rho(1-p)\geq 1$
we must have $z_\infty \geq \theta_2$ which implies
\[
z_\infty =\rho p\theta_1+\rho (1-p)\theta_2.
\]
 From the above we can conclude that
$P_s=0$ iff  $(1-\rho(1-p)) \theta_2< \rho p \theta 1$.
If the reverse inequality holds then all customers
of type 2 are being served successfully, i.e. $P_s=(1-p)$.
 \item {\em $\Theta$ exponentially distributed.} Assume w.l.o.g.\ that
the mean of $\Theta$ equals 1.
In this case $z_\infty $ can be determined from the equation
$z_\infty =\rho(1-{\rm e}^{-z_\infty})$ and $P_s={\rm e}^{-z_\infty}=1-z_\infty/\rho$.
\end{itemize}

Since $P_s$ does not depend on the mean of $\Theta$, and since the worst-case
property of the case of constant $\Theta$, it seems natural to conjecture
that the system performance is positively related to the variability
of $\Theta$. Thus it seems worthwhile to look
for ordering relations for $P_s$ if $\Theta_1 \stackrel{cvx}{\geq} \Theta_2$.
If $\expect{\Theta_1}=\expect{\Theta_2}$, this is
equivalent to $\expect{\min \{x,\Theta_1\}}\leq \expect{\min\{x,\Theta_2\}}$
for all $x\geq 0$.

Thus, if $\Theta_1 \stackrel{cvx}{\geq} \Theta_2$, it follows that
$z_{2,\infty}\geq z_{1,\infty}$ i.e.\ less variability in reneging behavior implies a lower service rate.
To prove that also $\pr{\Theta_1>z_{1,\infty}}\geq \pr{\Theta_2 > z_{2,\infty}}$ seems
hard without imposing further assumptions.

\subsection{Independent lead times}

In this case we can write (\ref{fixedpointeq2}) as
\[%\label{fpind}
\lambda \int_0^\infty \pr{B\geq u} \pr{D\geq z_\infty u} {\rm d} u = 1.
\]
which, in case $\expect{B}<\infty$, is equivalent to
$\pr{D\geq  z_\infty B^*}=1/\rho$,
with $B^*$ a random variable with density $\pr{B\geq x}/\expect{B}$.

Recall that  $P_s=\pr{D\geq z_\infty B}$.  Consequently, if $B$  is exponentially  distributed, we
have  the  insensitivity (w.r.t.\  the  distribution  of  $D$) result  $P_s=1/\rho$.   The
inequality  $P_s \leq  1/\rho$ holds  if  $B^*$ is  stochastically dominated  by $B$,  and
$P_s\geq 1/\rho$ vice versa.  Since $B^*$ being stochastically dominated by $B$ is related
to a low variability of $B$, we see  again that more variability (this time in the service
times) leads to a better system performance (i.e. higher $P_s$). \\

\noindent
{\em Exponential reneging}\\
If we assume that $D$ has an exponential distribution (and $B$ a general distribution), we see
that $z_\infty $ is the solution of
\beeq
\label{expreneg}
\rho \beta^*(z_\infty \nu))=1,
\eneq
with $\beta^*(s)=\expect{{\rm e}^{-sB^*}}$.

In addition, we have the following remarkable expression for the
complete fluid limit $z(t),t\geq 0$, if $z_0=0$:

\begin{prop}
Suppose $\pr{D\geq t}={\rm e}^{-\nu t}$, that $B$ is independent of $D$
and that $z_0=0$. Then the unique solution of (\ref{q}) is given by
\begin{equation}
\label{exponentialqueue}
z(t)=z_\infty (1-{\rm e}^{-\nu t}),
\end{equation}
with $z$ the solution of Equation (\ref{expreneg}).
\end{prop}

\begin{proof}
Recall that Equation (\ref{q}) has a unique solution.
We show that (\ref{exponentialqueue}) is indeed the solution of (\ref{q})
by verification.
We thus compute the right hand side of (\ref{q}) writing
$z(u)=z_\infty (1-{\rm e}^{-\nu u})$.

Observe that
\[
z_\infty \int_s^t \frac{1}{z(u)} \,{\rm d}u= \log ({\rm e}^{\nu t} -1)-\log ({\rm e}^{\nu s} -1).
\]
Consequently,
\begin{align*}
 \lambda \int_0^t \pr{D\geq t-s}& \pr{B\geq  \int_s^t (1/z(u)) {\rm d} u} \,{\rm d}s \\
&= \frac \lambda\nu {\rm e}^{-\nu t} \int_0^t\pr{z_\infty B \geq  \log ({\rm e}^{\nu t} -1)-\log ({\rm e}^{\nu s} -1)} d{\rm e}^{\nu s}\\
%&= \frac \lambda\nu {\rm e}^{-\nu t} \int_0^{{\rm e}^{\nu t} -1}\pr{z_\infty \nu B \geq  \log ({\rm e}^{\nu t}-1)-\log y}{\rm d}y \\
&= \frac \lambda\nu {\rm e}^{-\nu t}
\int_{-\log ({\rm e}^{\nu t}-1)}^\infty {\rm e}^{-v}
\pr{z_\infty \nu B \geq  \log ({\rm e}^{\nu t}-1) + z_\infty }{\rm d}v \\
&=
\frac \lambda\nu {\rm e}^{-\nu t} ({\rm e}^{\nu t}-1)
\int_0^\infty \pr{z_\infty \nu B\geq  v} {\rm e}^{-v} {\rm d}v\\
%&=   z_\infty (1-{\rm e}^{-\nu t})   \frac{\rho}{\expect{B}} \int_0^\infty\pr{B\geq  v/z_\infty \nu }{\rm e}^{-z_\infty\nu \frac{v}{z_\infty\nu}} {\rm d}(v/z_\infty\nu) \\
&=  z_\infty (1-{\rm e}^{-\nu t})  \rho  \beta^* (z_\infty \nu) 
=  z_\infty (1-{\rm e}^{-\nu t}).
\end{align*}
Which shows that $z_\infty (1-{\rm e}^{-\nu t})$ satisfies
 (\ref{q}).

\end{proof}

\subsection{TCP-friendly traffic}

Assume that there exist independent random variables $B_1$ and $D_1$ with finite means
such that
\begin{align*}
(B,D) &= (B_1,\infty) \hspace{1cm} \mbox{with probability } p, \\
&= (\infty, D_1) \hspace{1cm} \mbox{with probability } 1-p.
\end{align*}

When we view PS as a way of modeling TCP, this example models the integration
of elastic (TCP) traffic and TCP friendly UDP traffic; see Key {\em et al.}\ \cite{KMBK04} for a
related model. The latter type of traffic is using the system for a certain amount of time,
regardless of the level of congestion.

The fixed point equation (\ref{fixedpointeq2}) for $q$ specializes to
\[
z_\infty =\lambda p \expect{z_\infty B_1} +\lambda (1-p) \expect{D_1},
\]
Consequently, if the stability condition $\lambda p \expect{B_1}$ is satisfied,
we see that
\[
z_\infty =\frac{\lambda (1-p) \expect{D_1} }{1-\lambda p \expect{B_1}}.
\]

\section{Tightness}\label{s.tight} \setcounter{equation}{0} 

In this section we prove the first part of Theorem \ref{t.main}, that is, we
show that the sequence of processes $\{\Zbr\cd, r\in \R\}$ is tight in
$\mathbf{D}([0,\infty),\M)$. The main results in this section implying
this property are the compact containment Lemma \ref{l.compact}, and
an oscillation inequality in Lemma \ref{l.osc}. To prove these results,
a number of further lemmas are developed.  Section \ref{s.GC} derives
a Glivenko-Cantelli theorem for the stochastic primitives.  Section
\ref{s.dyn} introduces a fluid scaled version of the dynamic equation
for $\Zbr\cd$. The compact containment property is derived in Section
\ref{s.stability}. Section \ref{s.regularity} serves as a preparation
for the oscillation bound. In particular, it is shown that $\Zbr (t)$
charges arbitrarily small mass to thin $L$-shaped sets. The oscillation
bound is then shown in Section \ref{s.osc}.

Throughout this section, it is assumed that the assumptions of Section
\ref{scn:scaling} hold.

\subsection{A Glivenko-Cantelli theorem}\label{s.GC}
An important preliminary result is the following functional
Glivenko-Cantelli theorem for the stochastic primitives. It will
be convenient to consider them together as a single, measure valued
arrival process.  For $r\in\R$ and $t\ge s\ge 0$, define the fluid scaled
measure valued arrival process by
\[%\label{e.Ldef}
 \Lbr(t)=
   \frac{1}{r}\sum_{i=1}^{r\Er(t)}\delta_{(B^r_i,D^r_ir^{-1})},
\]
and define the fluid scaled increment 
\begin{equation}
\label{L-incr}
\Lbr(s,t)  = \Lbr(t)-\Lbr(s).
\end{equation}
Note that $\Lbr\cd$ is a random element of $\mathbf{D}([0,\infty),\M)$
and, for each $t\ge s\ge0$, $\Lbr(s,t)$ is a random element of $\M$.

To state and prove the result, we first introduce some notions from
empirical process theory. Our primary reference is \cite{VaaWel1996}.
A collection $\mathcal{C}$ of subsets of $\overline{\BR}_+^2$ {\em shatters}
an $n$-point subset $\{x_1,\dots,x_n\}\subset\overline{\BR}_+^2$ if the
collection $\{C\cap\{x_1,\dots,x_n\}:C\in\mathcal{C}\}$ has cardinality
$2^n$. In this case, say that $\mathcal{C}$ {\em picks out} all subsets
of $\{x_1,\dots,x_n\}$. The {\em Vapnik-\v{C}ervonenkis index (VC-index)}
of $\mathcal{C}$ is
\begin{equation}\nonumber
    V_\mathcal{C}=\min\{n:\text{$\mathcal{C}$ shatters no $n$-point
    subset}\},
\end{equation}
where the minimum of the empty set equals infinity.  The collection
$\mathcal{C}$ is a {\em Vapnik-\v{C}ervonenkis class (VC-class)}
if it has finite VC-index.

VC-classes satisfy a useful entropy bound. Let $\mathcal{Q}$ denote
the set of Borel probability measures on $\overline{\BR}_+^2$ and, for
$Q\in\mathcal{Q}$, let $\|f\|_Q=\langle |f|,Q\rangle$ denote the
$L_1(Q)$-norm of a Borel measurable function $f:\overline{\BR}_+^2\ar\BR$. For
$\varepsilon>0$, the $L_1(Q)$ $\varepsilon$-ball around $f$ is the
set of Borel functions $\{g:\|f-g\|_Q<\varepsilon\}$. For a family
of functions $\mathscr{V}$, the $(\varepsilon,L_1(Q))$-covering
number $N(\varepsilon,\mathscr{V},L_1(Q))$ is the smallest number
of $L_1(Q)$ $\varepsilon$-balls needed to cover $\mathscr{V}$.
If $\mathcal{C}$ is a VC-class, then for all $\varepsilon>0$, the family
$\mathscr{V}=\{1_C:C\in\mathcal{C}\}$ satisfies
\begin{equation}\label{e.finiteEntropy}
    \sup_{Q\in\mathcal{Q}}
    \log N(\varepsilon,\mathscr{V},L_1(Q))<\infty;
\end{equation}
see Theorem 2.6.4 in \cp{VaaWel1996}.

Recall the collection of corner sets $\mathcal{C}$ defined in Section
\ref{s.fluidModel}:
\begin{equation}\nonumber
    \mathcal{C}=\left\{ [x,\infty)\times[y,\infty):x,y\in\Rp \right\}
        \cup\left\{ [x,\infty]\times[y,\infty]:x,y\in\overline{\BR}_+ \right\}.
\end{equation}
Note that for any $3$-point subset $\{x_1,x_2,x_3\}\subset\overline{\BR}_+^2$,
it is impossible for $\mathcal{C}$ to pick out all three $2$-point subsets
of $\{x_1,x_2,x_3\}$. Since $\mathcal{C}$ shatters no $3$-point subset, it
has VC-index bounded above by $3$. Thus, $\mathcal{C}$ is a VC-class and
$\mathscr{V}=\{1_C:C\in\mathcal{C}\}$ satisfies \eqref{e.finiteEntropy}.

Define an envelope function for $\mathscr{V}$ as follows. Let
$\pi:\overline{\BR}_+^2\ar\overline{\BR}_+$ be the map $\pi(x,y)=\max\{x,y\}$. Since
$\pi$ is continuous, \eqref{a.lawConv} and the Skorohod representation
theorem imply the existence of $\overline{\BR}_+$-valued random variables
$X^r\sim\breve{\vartheta}^r\circ\pi^{-1}$ and $X\sim\vartheta\circ\pi^{-1}$
such that $X^r\ar X$ almost surely. Thus, there exists an $\overline{\BR}_+$-valued
random variable $Y$ such that 
\begin{equation}\label{e.Y}
	Y=\sup_{r\in\R}X^r, \quad\text{almost surely.}
\end{equation}
Let $\mu$ be the law of $Y$ on $\overline{\BR}_+$. Since $L_2(\mu)$
contains continuous unbounded functions, there exists a
continuous, unbounded function $\psi:\overline{\BR}_+\ar\BR_+$ that is
increasing on $[0,\infty)$, satisfies $\psi\ge 1$, and such that
$\langle\psi^2,\mu\rangle<\infty$. This implies that
\begin{equation}\label{e.Y2}
	\langle(\psi\circ\pi)^2,\vartheta\rangle
	=\BE[\psi(X)^2]
	\le \BE[\psi(Y)^2]
	< \infty.
\end{equation}
Let $F=\psi\circ\pi$, and note that $1_C\le F$ for all
$C\in\mathcal{C}$. That is, $F$ is an {\em envelope function} for
$\mathscr{V}$. Finally, define $\overline{\mathscr{V}}=\mathscr{V}\cup\{F\}$.

\begin{lem}\label{l.glivenko} 
Let $T>0$.  Then as $r\rightarrow\infty$,
\begin{equation}
\sup_{f\in\overline{\mathscr{V}}}\sup_{0\le s\le t\le T}
 \left|\langle f,\Lbr(s,t)\rangle
      - \lambda^r(t-s)\langle f,\breve{\vartheta}^r\rangle\right|
 \overset{\BP^r}{\lar} 0.
\end{equation}
\end{lem}

\begin{proof}
 Let $\varepsilon>0$. By \eqref{L-incr}, it suffices to show that
\[%\label{e.GC.suff}
    \limsup_{r\ar\infty}\BPr\left(
    \sup_{f\in\overline{\mathscr{V}}}\sup_{t\in[0,T]}
    \left|\langle f,\Lbr(t)\rangle
		-\lambda^rt\langle f,\breve{\vartheta}^r\rangle\right|
	>\varepsilon\right) \le \varepsilon.
\]
Note that the above event is measurable for each $r$ because it can
be rewritten using the suprema over rational $t$, and $f=1_C$ with $C$
having rational or infinite corner coordinates $x$ and $y$. Since $\langle
f,\Lbr(t)\rangle$ and $\lambda^rt\langle f,\breve{\vartheta}^r\rangle$
are nondecreasing in $t$ for each fixed $f\in\overline{\mathscr{V}}$, it
suffices to show that for each fixed $t\in[0,T]$,
\begin{equation}
    \limsup_{r\ar\infty}\BPr\left( \sup_{f\in\overline{\mathscr{V}}}
        \left|\langle f,\Lbr(t)\rangle
			-\lambda^rt\langle f,\breve{\vartheta}^r\rangle\right|
        >\varepsilon\right) \le\varepsilon.
\nonumber
\end{equation}
Since
\begin{equation}
 \langle f,\Lbr(t)\rangle-\lambda^rt\langle f,\breve{\vartheta}^r\rangle
  =\langle f,\breve{\vartheta}^r\rangle\left(\Er(t)-\lambda^rt\right)
    + \Er(t)\left(\frac{\langle f,\Lbr(t)\rangle}{\Er(t)}
        -\langle f,\breve{\vartheta}^r\rangle\right)
\nonumber
\end{equation}
(with the convention that division by zero equals zero), it suffices to
show the two bounds
\begin{align}%\label{e.GC.suffA}
  \limsup_{r\ar\infty}\BPr\left(\sup_{f\in\overline{\mathscr{V}}}
    \left|\langle f,\breve{\vartheta}^r\rangle
		\left(\Er(t)-\lambda^rt\right)\right|
    >\frac{\varepsilon}{2}\right)\le\frac{\varepsilon}{2},\notag \\
    \limsup_{r\ar\infty}\BPr\left(\sup_{f\in\overline{\mathscr{V}}}
    \left|\Er(t)\left(\frac{\langle f,\Lbr(t)\rangle}{\Er(t)}
            -\langle f,\breve{\vartheta}^r\rangle\right) \right|
        >\frac{\varepsilon}{2}\right)\le\frac{\varepsilon}{2}.
\label{e.GC.suffB}
\end{align}
The first equation follows from assumption~\eqref{a.u-conv-exp} and by observing that 
\begin{equation}
	\sup_{r\in\R}\sup_{f\in\overline{\mathscr{V}}}
		\langle f,\breve{\vartheta}^r\rangle
	\le \sup_{r\in\R}\langle F,\breve{\vartheta}^r\rangle
	= \sup_{r\in\R}\BE[\psi(X^r)]
	\le \BE[\psi(Y)]<\infty,
\label{e.GC.envelopeBound}
\end{equation}
which follows from \eqref{e.Y} and \eqref{e.Y2}.  To show
\eqref{e.GC.suffB}, it suffices to verify three assumptions of
Theorem 2.8.1 in \cp{VaaWel1996}.  Observe that for each $n\in\BN$
and $(e_1,\dots,e_n)\in\BR^n$, the function
\begin{equation}
\nonumber
 (x_1,\dots,x_n)\to\sup_{f\in\overline{\mathscr{V}}}
	\sum_{i=1}^ne_if(x_i)
\end{equation}
is measurable on the completion of $(\overline{\BR}_+^2, \mathscr{B},
\breve{\vartheta}^r)^n$, for each $r\in\R$.  Thus, $\overline{\mathscr{V}}$
is a {\em $\breve{\vartheta}^r$-measurable class} for each $r\in\R$;
see Definition 2.3.3 in \cp{VaaWel1996}.  Moreover, $\overline{\mathscr{V}}$
is uniformly bounded above by the envelope function $F$, and 
\begin{equation}
	\lim_{M\ar\infty}\sup_{r\in\R}
 	\langle F1_{\{F>M\}},\breve{\vartheta}^r\rangle=0,
\nonumber
\end{equation}
by Markov's inequality, \eqref{e.Y}, and \eqref{e.Y2}.
Lastly, $\overline{\mathscr{V}}$ satisfies the finite entropy bound
\eqref{e.finiteEntropy} because $N(\varepsilon,\overline{\mathscr{V}},L_1(Q))
\le N(\varepsilon,\mathscr{V},L_1(Q))+1$ and $\mathcal{C}$ is a
VC-class. The previous three observations imply that the assumptions
of Theorem 2.8.1 in \cp{VaaWel1996} are satisfied. Consequently,
$\overline{\mathscr{V}}$ is {\em Glivenko-Cantelli, uniformly in $r$}. That
is, for every $\delta>0$, there exists an $n_\delta$ such that $n\ge
n_\delta$ implies
\begin{equation}
    \sup_{r\in\R}\BP^r\left(
    \sup_{m\ge n}\sup_{f\in\overline{\mathscr{V}}}
    \frac{1}{m}\sum_{i=1}^m f(B^r_i,D^r_ir^{-1})
       -\langle f,\breve{\vartheta}^r\rangle
    >\delta \right)  \le \delta.
\label{e.GV.uniformGV}
\end{equation}
Choose $\delta=\min\{\varepsilon/2,\varepsilon/(4\lambda T)\}$. 
The left side of \eqref{e.GC.suffB} is bounded above by
\begin{equation}
\nonumber
\limsup_{r\ar\infty}\BPr\left(\Er(t)>2\lambda T\right)
  +\limsup_{r\ar\infty}\BPr\left(\sup_{f\in\overline{\mathscr{V}}}
    \left|\frac{\langle f,\Lbr(t)\rangle}{\Er(t)}
            -\langle f,\breve{\vartheta}^r\rangle \right|
        >\frac{\varepsilon}{4\lambda T}\right).
\end{equation}
The first term equals zero by \eqref{a.u-conv-exp}.  For the second
term, rewrite
\begin{equation}
    \frac{\langle f,\Lbr(t)\rangle}{\Er(t)}=\frac{1}{E^r(rt)}
       \sum_{i=1}^{E^r(rt)}f(B^r_i,D^r_ir^{-1}),
\nonumber
\end{equation}
and bound each probability in the second term by 
\begin{equation}
\label{e.GC.almost}
 \BP^r(E^r(rt)< n_{\delta})
+ \BP^r\left(
    \sup_{m\ge n_\delta}\sup_{f\in\overline{\mathscr{V}}}
    \frac{1}{m}\sum_{i=1}^mf(B^r_i,D^r_ir^{-1})
		-\langle f,\breve{\vartheta}^r\rangle
    >\frac{\varepsilon}{4\lambda T} \right).
\end{equation}
By \eqref{a.u-conv-exp}, the first term in \eqref{e.GC.almost} converges
to zero as $r\rightarrow\infty$.  By \eqref{e.GV.uniformGV}, the second
term is bounded above by $\delta\le\varepsilon/2$, uniformly in $r\in\R$.
This implies \eqref{e.GC.suffB}.  
\end{proof}

\subsection{Fluid scaled dynamic equation}
\label{s.dyn}

Using \eqref{e.dynamicEQ}, it is easy to see that the fluid scaled state
descriptor of the $r$th model satisfies the following equation almost
surely: for each Borel set $A\in\mathscr{B}$, and all $t,h\ge 0$,
\begin{multline}
\label{e.dyn}
  \Zbr(t+h)(A) =
   \Zbr(t)\left(A+(\Sbr (t,t+h),h)\right)
   \\+ \sum_{i=r\Er(t)+1}^{r\Er(t+h)}
      1^+_A\left(\vbr_i(t+h), \lbr_i(t+h)\right).
\end{multline}

Subsequent proofs use estimates obtained from this equation. Two
estimates result from bounding the summands in \eqref{e.dyn} by $1$
and optionally bounding the first term on the right side by its total
mass; for each $A\in\mathscr{B}$ and $t,h\ge 0$,
\begin{align}\label{e.dynm.bound1}
 \begin{split}
  \Zbr(t+h)(A) &\le\Zbr(t)\left(A+(\Sbr (t,t+h),h)\right)+\Lbr(t,t+h)(\overline\BR^2_+) \\
    & \le \Zbr(t)\left(\overline\BR^2_+\right)+\Lbr(t,t+h)(\overline\BR^2_+).
 \end{split}
\end{align}
Two more estimates follow from \eqref{e.dyn} by simply ignoring any
arrivals; for each $A\in\mathscr{B}$ and $t,h\ge 0$,
\begin{equation}\label{e.dynm.bound2}
 \Zbr(t)\paren{A+(\Sbr(t,t+h),h)}\le\Zbr(t+h)(A)\le\Zbr(t+h)\left(\overline \BR_+^2\right).
\end{equation}

\subsection{Compact containment}
\label{s.stability}

This section establishes the compact containment property needed to prove
tightness. 

\begin{lem}
\label{l.compact}
Let $T>0$ and $\eta>0$. There exists a compact set $\K\subset\M$
such that
\begin{equation}\label{l.compact.eq}
    \liminf_{r\ar\infty}\BPr\left( \Zbr(t)\in\K
      \,\,\text{for all $t\in[0,T]$} \right)\ge 1-\eta.
\end{equation}
\end{lem}

\begin{proof}
A set  $\K \subset\M$ is relatively compact if $\sup_{\xi \in \K }
\xi(\overline{\BR}_+^2 )<\infty$, and if there exists a sequence of nested
compact sets $K_n\subset\overline \BR_+^2$ such that $\bigcup_{n\in\BN}
K_n=\overline{\BR}_+^2$ and
$$
\lim_{n\rightarrow\infty} \sup_{\xi \in \K} \xi(K_n^c)=0,
$$
where $K_n^c$ denotes the complement of $K_n$; see \cite{Ka1976},
Theorem A 7.5.  Consider the nested sequence of compact sets in $\overline
\BR^2_+$ given by
$$
K_n = ([0,n]\times [0,n]) \cup ([0,n] \times \{\infty\}) 
	\cup (\{\infty\} \times [0,n]) 
	\cup (\{\infty\} \times \{\infty\}),\quad n\in\BN.
$$

By \eqref{a.initconv}, $\Zbr(0)\wk\zeta_0$ in distribution, and so
the sequence $\{\Zbr(0)\}$ is tight. Thus, there is a compact set
$\K_0\subset\M$, such that
\begin{equation}
\label{comp.z0}
\liminf_{r\rightarrow\infty} \BPr( \Zbr(0) \in \K_0 ) \geq
1-\frac{\eta}{2}.
\end{equation}
Let $M_0=\sup_{\xi\in\K_0}\xi(\overline{\BR}_+^2)$, and let
$a_n=\sup_{\xi\in\K_0}\xi(K_n^c)$ for each $n\in\BN$.  Since $\K_0$ is
compact, $M_0<\infty$ and there exists a sequence of nested compact sets
$J_n\subset\overline\BR_+^2$ such that $\bigcup_{n\in\BN}J_n=\overline{\BR}_+^2$
and $\lim_{n\ar\infty}\sup_{\xi\in\K_0} \xi(J_n^c)=0$.  Since $J_n\subset
K_{k(n)}$ for each $n\in\BN$ and sufficiently large $k(n)\in\BN$, it
follows that $a_n\ar0$ as $n\ar\infty$.

Recall the definition from Section \ref{s.GC} of the envelope
function $F=\psi\circ\pi$ for the family $\overline{\mathscr{V}}$.
By \eqref{a.u-conv-exp} and \eqref{e.GC.envelopeBound}, the constant
$M=\sup_{r\in\R}\left(\lambda^rT\langle F,\breve{\vartheta}^r\rangle+1
\right)$ is finite. 
Let $\K$ be the closure of the set
\begin{equation}
	\left\{\xi\in\M:\xi(\overline{\BR}_+^2)\le M_0+M
		\,\,\text{and}\,\,
	\xi(K_n^c)\le a_n+\psi(n)^{-1}M
				\,\,\text{for all $n\in\BN$}\right\}.
\nonumber
\end{equation}
Since $a_n+\psi(n)^{-1}M\ar0$ as $n\ar\infty$, the set $\K$ is compact
in $\M$.

For each $r\in\R$, denote the event in \eqref{comp.z0} by $\Omega^r_0$
and define the event
\[%\label{e.CC.Omega1n}
	\Omega^r_1=\left\{
	 \langle F,\Lbr(T)\rangle
		\le \lambda^rT\langle F,\breve{\vartheta}^r\rangle+1\right\}.
\]
By \eqref{comp.z0} and Lemma \ref{l.glivenko},
$\liminf_{r\ar\infty} \BPr(\Omega^r_0\cap\Omega^r_1)\ge 1-\eta$. Fix
$\omega\in\Omega^r_0\cap\Omega^r_1$ and $t\in[0,T]$, and assume for the
remainder of the proof that all random objects are evaluated at this
$\omega$. Then it suffices to show that $\Zbr(t)\in\K$.

By \eqref{e.dynm.bound1},
\begin{equation}
	\Zbr(t)(\overline{\BR}_+^2)\le \Zbr(0)(\overline{\BR}_+^2)
		+\Lbr(t)(\overline{\BR}_+^2).
\nonumber
\end{equation}
Since $\Lbr(t)(\overline{\BR}_+^2)=\langle 1,\Lbr(t)\rangle\le\langle
1,\Lbr(T)\rangle\le\langle F,\Lbr(T)\rangle$, the definitions of
$\Omega^r_0$, $\Omega^r_1$, and $M$ imply that
\begin{equation}
	\Zbr(t)(\overline{\BR}_+^2)\le M_0
		+M.
\label{e.CC.mass}
\end{equation}
Fix $n\in\BN$. By \eqref{e.dyn}, 
\begin{equation}
\nonumber
  \Zbr(t)(K_n^c) =
   \Zbr(0)\left(K_n^c+(\Sbr (0,t),t)\right)
   + \sum_{i=1}^{r\Er(t)}
      1^+_{K_n^c}\left(\vbr_i(t), \lbr_i(t)\right).
\end{equation}
The shape of the set $K_n^c$ implies that 
\[
K_n^c+(S(0,t),t)\subset
K_n^c \text{ and } 1^+_{K_n^c} \left(\vbr_i(t), \lbr_i(t)\right)\le
1_{K_n^c}(B^r_i,D^r_ir^{-1}),
\]
for $i=1,\dots,r\Er(t)$. Thus,
\begin{equation}
\nonumber
  \Zbr(t)(K_n^c) \le
   \Zbr(0)\left(K_n^c\right)
   + \langle 1_{K_n^c},\Lbr(t)\rangle.
\end{equation}
By definition of $\psi$, $F$, and by Markov's inequality, $1_{K_n^c}\le
\psi(n)^{-1}F$. So
\[
  \Zbr(t)(K_n^c) \le
   \Zbr(0)\left(K_n^c\right)
   + \psi(n)^{-1}\langle F,\Lbr(t)\rangle.
\]
Since $\langle F,\Lbr(t)\rangle\le\langle F,\Lbr(T)\rangle$, the
definitions of $\Omega^r_0$, $\Omega^r_1$, and $M$ imply that
\begin{equation}
\label{e.CC.tails}
  \Zbr(t)(K_n^c) \le
   a_n
   + \psi(n)^{-1}M.
\end{equation}
Equations \eqref{e.CC.mass} and \eqref{e.CC.tails} imply that
$\Zbr(t)\in\K$.
\end{proof}

\subsection{Asymptotic regularity}\label{s.regularity}

The second and main step necessary to prove tightness is to bound
the probability that the process $\Zbr\cd$ oscillates. Oscillations
may result from sudden arrivals or departures of a large amount of
mass. Sudden arrivals are controlled by the regularity of the arrival
process. To show that sudden departures are unlikely as well, we show that
$\Zbr\cd$ assigns arbitrarily small mass to the boundaries of the sets
$C\in\mathcal{C}$.  This is phrased in terms of $\kappa$-enlargements
of the boundaries of these sets (forming a collection of $L$-shaped
sets).  For $C\in\mathcal{C}$ and $\kappa>0$, let $\partial_C$ denote
the boundary of $C$ and let $$\partial^\kappa_C=\left\{w\in\overline
\BR^2_+:\inf_{z\in \partial_C}\|w-z\|<\kappa\right\}$$ be the
$\kappa$-enlargement in $\overline\BR^2_+$ of its boundary, where the
infimum over the empty set equals $\infty$. (Note that $\partial_C$,
and therefore also $\partial_C^\kappa$, is empty for the corner
sets $\overline{\BR}_+^2$ and $\{\infty\}\times\{\infty\}$. Note also that
$\partial_C^\kappa=((x-\kappa)^+,x+\kappa)\times\{\infty\}$ for a corner
set of the form $[x,\infty]\times\{\infty\}$ with $x\in[0,\infty)$.) The
following lemma establishes the result for the initial condition
$\Zbr(0)$.

\begin{lem}
\label{l.noatoms0}
For all $\varepsilon,\eta>0$ there exists a $\kappa>0$ such that
\begin{equation}
  \liminf_{r\ar\infty}\BPr\left(
\sup_{C\in\mathcal{C}}\Zbr(0)(\partial^\kappa_C)
  \le\varepsilon \right) \ge 1-\eta.
\label{e.noatoms0.claim}
\end{equation}
 \end{lem}
\begin{proof}
Fix $\varepsilon,\eta>0$ and let $\Zbr_1(0)(\cdot)=\Zbr(0)(\cdot \times
\overline\BR_+)$ and $\Zbr_2(0)(\cdot)=\Zbr(0)(\overline\BR_+ \times \cdot)$.
For each $C\in \mathcal{C}$ and $\kappa>0$, $$\partial_C^\kappa\subset
([x,x+2\kappa]\times\overline{\BR}_+)\cup(\overline{\BR}_+\times [y,y+2\kappa]),$$
for some $(x,y)\in\BR_+^2=[0,\infty)\times[0,\infty)$. Thus, it suffices
to show that, for $i=1,2$, there exists a $\kappa>0$ such that
\begin{equation}
\label{eq.noatoms0.goal}
  \liminf_{r\ar\infty}\BPr\left( 
	\sup_{x\in[0,\infty)}\Zbr_i(0)([x,x+2\kappa])
  \le\varepsilon \right) \ge 1-\frac{\eta}{2}.
\end{equation}
We prove the statement for $i=1$; the proof is identical
for $i=2$.  

The projection $(x,y)\mapsto x$ is continuous, so \eqref{a.initconv}
implies that $\Zbr_1(0)$ converges in distribution to
$\zeta_0(\cdot\times\overline{\BR}_+)$ as $r\ar\infty$.  Since
$\zeta_0(\cdot\times\overline{\BR}_+)$ is free of atoms in $[0,\infty)$,
there exists a $\kappa>0$ such that
\begin{equation}
\label{e.noatoms0.grid}
 \sup_{x\in [0,\infty)}\zeta_0([x,x+4\kappa]\times\overline{\BR}_+)
\le \varepsilon. 
\end{equation}
(If \eqref{e.noatoms0.grid} fails, it is easy to construct an atom of
$\zeta_0(\cdot\times\overline{\BR}_+)$.)  Moreover, there exists a constant
$M$ such that
\begin{equation}
\label{eq.noatoms.tail}
 \zeta_0([M,\infty)\times\overline{\BR}_+) \le\varepsilon.
\end{equation}

Let $N=\lceil M/\kappa\rceil+1$, where $\lceil x\rceil$ denotes
the smallest integer $n\ge x$. For $n=1,\ldots, N-1$, define  the
set $I_n=[n\kappa, (n+4)\kappa]$ and define $I_N=[M,\infty)$.
Note that, for every $x\in[0,\infty)$ there is an $n\le N$ such that
$[x,x+2\kappa]\subset I_n$.  To prove \eqref{eq.noatoms0.goal}, it
therefore suffices to show that
\begin{equation}
\label{eq.noatoms0.goal2}
  \liminf_{r\ar\infty} \BPr\left( \max_{n\le N}\Zbr_1(0)(I_n)
  \le\varepsilon \right) \ge 1-\frac{\eta}{2}.
\end{equation}

Let $\M(\overline{\BR}_+)$ denote the space of finite nonnegative Borel
measures on $\overline{\BR}_+$, endowed with the weak topology. Let
$\mathbf{A}=\left\{ \xi\in\M(\overline{\BR}_+): \max_{n\le
N}\xi(I_n)<\varepsilon \right\}$, and suppose that a sequence
$\{\xi_k\}\subset\M(\overline{\BR}_+)$ satisfies $\xi_k\wk\xi$ for some
$\xi\in\mathbf{A}$. Since the sets $I_n$ are closed, the Portmanteau
theorem (adapted to finite measures) implies that
\begin{equation}
\limsup_{k\rightarrow\infty} \xi_k(I_n) \le \xi(I_n)<\varepsilon,
	\quad\text{for all $n\le N$}. 
\nonumber
\end{equation}
Hence, $\xi_k \in \mathbf{A}$ for sufficiently large $k$, which implies
that $\mathbf{A}$ is open in $\M(\overline{\BR}_+)$. Thus, a second application
of the Portmanteau theorem yields
\begin{equation}
\nonumber
\liminf_{r\rightarrow\infty} \BPr(\Zbr_1(0) \in \mathbf{A}) \ge
\BP(\zeta_0(\cdot\times\overline{\BR}_+) \in \mathbf{A}) = 1,
\end{equation}
which implies \eqref{eq.noatoms0.goal2}. \end{proof}

The regularity result is now shown for the entire state descriptor
$\Zbr\cd$.  

\begin{lem}
\label{l.noatoms}
Let $T>0$ and $\varepsilon,\eta>0$. There exists a $\kappa>0$ such that
\begin{equation}
\label{e.noatoms}
 \liminf_{r\ar\infty}\BPr\left(
   \sup_{C\in\mathcal{C}}\sup_{t\in [0,T]}
   \Zbr(t)(\partial_C^\kappa) \le \varepsilon  \right)
      \ge 1-\eta.
\end{equation}
\end{lem}

\begin{proof}
 By  Lemmas \ref{l.glivenko}, \ref{l.compact}, and \ref{l.noatoms0},
there exists a compact $\K\subset\M$ and a $\kappa_0>0$,  such that for
all $\delta>0$, the events
\begin{align}
\notag
\begin{split}
 \Omega^r_1 & =\left\{\sup_{C\in\mathcal{C}}
        \Zbr(0)(\partial_C^{\kappa_0})
             \le\frac{\varepsilon}{2}\right\},\\
  \Omega^r_2  &= \left\{\sup_{C\in\mathcal{C}}\sup_{0\le s\le t\le T}
        \left|\Lbr(s,t)(C) - \lambda^r(t-s)\breve\vartheta^r (C)\right|
                \le \delta \right\},\\
  \Omega^r_3 & = \left\{\Zbr(t)\in\K\,\,
      \text{for all $t\in [0,T]$} \right\},\\
 \Omega^r_0
    & =\Omega^r_1\cap \Omega^r_2\cap\Omega^r_3,
\end{split}
\end{align}
satisfy
\begin{equation}
\label{e.noatoms.goodEvent}
 \liminf_{r\ar\infty}\BPr\paren{\Omega^r_0}\ge 1-\eta.
\end{equation}
Recall the compact sets $K_n$ defined in the proof of Lemma
\ref{l.compact}.  Since $\K$ is compact, there exists a finite $M\ge 1$
and an integer $R<\infty$ such that
\begin{align}
\label{e.noatoms.massBound}
    \sup_{\xi\in\K} \xi(\overline{\BR}_+^2) &\le M, \\
\label{e.noatoms.noTail}
    \sup_{\xi\in\K}\xi\left(K_R^c\right) &\le
    \frac{\varepsilon}{2}.
\end{align}
Let $\lambda^*=\sup_{r\in\R}\lambda^r$, which is finite by
\eqref{a.u-conv-exp}.  Fix 
\[
h=\varepsilon(8\lambda^*)^{-1}, \quad
\kappa=\min\{\kappa_0,h(2M)^{-1}\} \text{  and }\delta=\varepsilon
\min\{(8\lceil RMh^{-1}\rceil)^{-1}, 2^{-1}\}.
\]
  For $r\in\R$,
let $\Omega^r_*$ denote the event in \eqref{e.noatoms}.
By \eqref{e.noatoms.goodEvent}, it suffices to show that
$\Omega^r_0\subset\Omega^r_*$.  Let $\omega\in\Omega^r_0$ be arbitrary;
for the remainder of the proof, all random objects are evaluated at
this $\omega$.

Consider any $r\in\R$, $t\in [0,T]$ and $C\in\mathcal{C}$. We must show
that $\Zbr(t)(\partial_C^\kappa)\le\varepsilon$.  Define the random time
\begin{equation}\nonumber
 \tau_1  =\sup\{s\le t:\langle 1,\Zbr(s)\rangle=0\},
\end{equation}
if the supremum exists, and define $\tau_1=0$ otherwise. Let
$\tau=\max\{\tau_1,t-RM\}$. We first show that
\begin{equation}\label{e.noatoms.tau}
  \Zbr(\tau)\left( \partial_C^\kappa+(\Sbr (\tau,t),t-\tau)\right)
     \le \frac{\varepsilon}{2}.
\end{equation}
If $\tau=0$, this follows from the definition of $\Omega^r_1$
because $\kappa\le\kappa_0$, because
\begin{equation}
\notag
\partial_C^\kappa+(\Sbr (\tau,t),t-\tau) \subset
\partial_{C+(\Sbr (\tau,t),t-\tau)}^\kappa,
\end{equation}
and because $\mathcal{C}$ is closed under positive translation.
Suppose $\tau=\tau_1>0$. Then there is a sequence $\{\tau_n\}$, with
$\tau_n\uparrow\tau$, such that $\langle 1,\Zbr(\tau_n)\rangle=0$ for
all $n$.  In this case, \eqref{e.dynm.bound1} and the definition of
$\Omega^r_2$ imply that, for all $n$,  
\begin{equation}
\nonumber
 \Zbr(\tau)\left( \partial_C^\kappa+(\Sbr (\tau,t),t-\tau)\right)\le
    \Zbr(\tau_n)\left(\overline \BR_2^+\right)
		+ \Lbr(\tau_n,\tau)\left(\overline \BR_2^+\right)
    \le \lambda^r(\tau-\tau_n) + \delta.
\end{equation}
Letting $\tau_n\uparrow \tau$ yields
\begin{equation}
\nonumber
 \Zbr(\tau)\left( \partial_C^\kappa+(\Sbr (\tau,t),t-\tau)\right)\le
  \delta \le \f{\varepsilon}{2}.
\end{equation}
Suppose that $\tau=t-RM$.  Since $\langle 1,\Zbr(s)\rangle>0$
for all $s\in(\tau,t]$, the definition of $\Omega^r_3$ and
\eqref{e.noatoms.massBound} imply that
\begin{equation}
\notag
  \Sbr (\tau,t)=\int_{t-RM}^t\langle 1,\Zbr(s)\rangle^{-1} \,{\rm d}s
    \ge R.
\end{equation}
Thus, by the definition of $\Omega^r_3$ and \eqref{e.noatoms.noTail},
\begin{equation}
\notag
    \Zbr(\tau)\left( \partial_C^\kappa+(\Sbr (\tau,t),t-\tau)\right)
     \le \Zbr(\tau)\left(K_R^c\right)\le \frac{\varepsilon}{2},
\end{equation}
which proves \eqref{e.noatoms.tau}.

By \eqref{e.dyn}, 
\begin{multline}\label{e.noatoms.1}
 \Zbr(t)\left( \partial_C^\kappa\right) = \Zbr(\tau)
   \left( \partial_C^\kappa+(\Sbr (\tau,t),t-\tau)\right)
\\    +\f{1}{r}\sum_{i=r\Er(\tau)+1}^{r\Er(t)}
     1^+_{\partial_C^\kappa}
    \left(\overline B^r_i(t),\overline D^r_i(t)\right).
\end{multline}
Let $I$ denote the second right hand term in \eqref{e.noatoms.1}.
By \eqref{e.noatoms.tau}, it remains to show that $I\le \varepsilon/2$.
Let $N=\lceil(t-\tau)h^{-1}\rceil$ and, for each $n=0,\ldots,N-1$,
let $t_n=\tau+nh$ and $t^n=\min\{t_{n+1},t\}$. Then, using
the inequality $1^+_{\partial_C^\kappa}(\cdot,\cdot)\le
1_{\partial_C^\kappa}(\cdot,\cdot)$,
\begin{equation}
\label{e.noatoms.firstI}
    I\le\sum_{n=0}^{N-1}\frac{1}{r}\sum_{i=r\Er(t_n)+1}^{r\Er(t^n)}
       1_{\partial_C^\kappa}
         \left(\overline B^r_i(t),\overline D^r_i(t)\right).
\end{equation}
Consider $n\in\{0,\ldots,N-1\}$ and $i$ such that $U^r_ir^{-1}\in
(t_n,t^n]$. Observe that
\begin{equation}
\label{e.noatoms.Sbound}
    \Sbr (t^n,t)
        \le\Sbr (U^r_ir^{-1},t)
          \le\Sbr (t_n,t).
\end{equation}
By definition,
\begin{equation}
\label{e.noatoms.1rewrite}
    1_{\partial_C^\kappa}(\overline B^r_i(t),\overline D^r_i(t))
     = 1_{\partial_C^\kappa
          +(\Sbr (U^r_ir^{-1},t),t-U^r_ir^{-1})}(B^r_i,D^r_ir^{-1}).
\end{equation}
So, letting
\begin{equation}
\notag
\begin{split}
    C_n^- &=C+\left(\Sbr (t^n,t)-\kappa,
         t-t^n-\kappa\right) \cap \overline \BR_2^+, \\
    C_n^+ &=C+\left(\Sbr (t_n,t)+\kappa,
         t-t_n+\kappa\right)\cap \overline \BR_2^+,\\
    C_n &=C_n^-\setminus C_n^+,
\end{split}
\end{equation}
it follows from \eqref{e.noatoms.Sbound} and \eqref{e.noatoms.1rewrite}
that
\begin{equation}
\label{e.noatoms.1bound}
    1_{\partial_C^\kappa}(\overline B^r_i(t),\lbr_i(t))
      \le 1_{C_n}(B^r_i,D^r_ir^{-1}).
\end{equation}
Conclude from \eqref{e.noatoms.firstI} and \eqref{e.noatoms.1bound} that
\begin{multline*}
 I\le\sum_{n=0}^{N-1}\frac{1}{r}
     \sum_{i=r\Er(t_n)+1}^{r\Er(t^n)}1_{C_n}(B^r_i,D^r_ir^{-1}) 
	\\=\sum_{n=0}^{N-1}\left(
       \Lbr (t_n,t^n)(C_n^-)-\Lbr (t_n,t^n) (C_n^+)\right).
\end{multline*}
For all $n<N$, $C_n^-,C_n^+\in\mathcal{C}$ and $t^n-t_n\le h$.  So the
definition of $\Omega^r_2$ implies that
\begin{equation}
\notag
    I\le\sum_{n=0}^{N-1}\left(\lambda^r h \breve{\vartheta}^r(C_n) 
		+2\delta\right).
\end{equation}
By definition of $N$, and since $t-\tau\le RM$,
\begin{equation}
\notag
    I\le\lambda^* h
        \sum_{n=0}^{N-1} \breve{\vartheta}^r(C_n)
           + \lceil RMh^{-1}\rceil 2\delta.
\end{equation}
This implies, by choice of $\delta$, that
\begin{equation}
\label{e.noatoms.evenOdd}
    I\le\lambda^* h
        \sum_{n=0}^{N-1} \breve{\vartheta}^r(C_n)
           +\frac{\varepsilon}{4}.
\end{equation}
If $n\in\{0,\ldots,N-3\}$, then
\begin{equation}
\Sbr (t_{n+1},t_{n+2})\ge hM^{-1}\ge 2\kappa,
\nonumber
\end{equation}
because $0<\langle 1,\Zbr(s)\rangle\le M$ for all $s\in (\tau,t]$
and because $h\ge \kappa 2M$ by definition. Thus, for all
$n\in\{0,\ldots,N-3\}$,
\begin{equation}
\nonumber
     \Sbr (t^n,t)-\kappa =\Sbr (t_{n+1},t_{n+2})
                    +\Sbr (t_{n+2},t)-\kappa
     \ge \Sbr (t_{n+2},t)+\kappa.
\end{equation}
Hence, $C_n^-\subset C_{n+2}^+$ for all $n\in\{0,\ldots,N-3\}$,
and consequently, $C_n\cap C_{n+2}=\varnothing$. Thus,
since $\breve{\vartheta}^r$ is a probability measure,
\[
\sum_{n=0}^{\lfloor (N-1)/2\rfloor}\breve{\vartheta}^r(C_{2n}) \text{ and }
\sum_{n=0}^{\lfloor(N-2)/2\rfloor}\breve{\vartheta}^r(C_{2n+1})
\]
are both bounded above by one. Conclude from \eqref{e.noatoms.evenOdd} that
\begin{equation}
\notag
    I\le 2\lambda^* h+\frac{\varepsilon}{4},
\end{equation}
which implies, by choice of $h$, that $I\le\varepsilon/2$.
\end{proof}

\subsection{Oscillation bound}\label{s.osc}

This section establishes the second main ingredient for proving tightness
of the state descriptors.  As a metric on $\M$, we use the Prohorov metric
(adapted to finite measures). For $\mu,\nu\in\M$, define
\[
\dw{\mu,\nu}{=}
  \inf\left\{\rule{0mm}{5mm}\varepsilon{>}0: 
  \mu(A){\le} \nu(A^\varepsilon)+\varepsilon \text{ and }
  \nu(A){\le} \mu(A^\varepsilon)+\varepsilon\,\,\,
 	\text{for all closed $A\in \mathscr{B}$}\right\}.
\]
Recall that $A^\varepsilon=\{w\in\overline{\BR}_+^2:\inf_{z\in
A}\|z-w\|<\varepsilon\}$ and that $\mathscr{B}$ denotes the Borel subsets
of $\overline{\BR}_+^2$.

\begin{defn}
For each $\zeta\cd\in \mathbf{D}([0,\infty),\M)$ and each $T>\delta>0$,
define the modulus of continuity on $[0,T]$ by
\[%\label{e.modc}
\ww{\zeta\cd,\delta}= \sup_{t\in [0,T-\delta]}
\sup_{h\in[0,\delta]} \dw{\zeta(t+h),\zeta(t)}.
\]
\end{defn}

\begin{lem}\label{l.osc}
For all $T>0$ and $\varepsilon,\eta\in(0,1)$, there exists
$\delta\in(0,T)$ such that
\begin{equation}\label{e.osc}
\liminf_{r\ar\infty}
  \BPr \paren{ \ww{\Zbr\cd,\delta}\le\varepsilon}
      \ge 1-\eta.
\end{equation}
\end{lem}

\begin{proof}
 As before, let $\lambda^*=\sup_{r\in\R}\lambda^r.$
For each $\kappa>0$, define
\begin{equation}
\nonumber
L_\kappa=([0,\kappa]\times\overline{\BR}_+)\cup(\overline{\BR}_+\times[0,\kappa]).
\end{equation}
By Lemmas \ref{l.glivenko} and \ref{l.noatoms}, there exists
$\kappa\in(0,1)$ such that for all $\delta\in(0,T)$, the events
\begin{align}
\nonumber
 \begin{split}
  \Omega^r_1 &=
         \left\{\sup_{t\in [0,T]}
         \Zbr(t)(L_\kappa)  \le\frac{\varepsilon}{4}\right\}, \\
  \Omega^r_2 & =
         \left\{\sup_{t\in [0,T-\delta]}
        \Lbr (t,t+\delta)(\overline \BR_+^2)\le 2\lambda^*\delta\right\}, \\
  \Omega^r_0&= \Omega^r_1\cap\Omega^r_2,
 \end{split}
\end{align}
satisfy
\begin{equation}\label{e.osc.event}
 \liminf_{r\ar\infty}\BPr\left(\Omega^r_0\right)\ge 1-\eta.
\end{equation}
Fix $\delta =\kappa\varepsilon^2(8\max\{\lambda^*,1\})^{-1}$ and let
$\Omega^r_*$ be the event in \eqref{e.osc}. By \eqref{e.osc.event}, it
suffices to show that $\Omega^r_0\subset \Omega^r_*$ for each $r$.  Fix
$r\in\R$ and $\omega\in \Omega^r_0$; for the remainder of the proof all
random objects are evaluated at this $\omega$.  Fix  $t\in[0,T-\delta]$,
$h\in[0,\delta]$ and let $A\in\mathscr{B}$ be closed. It suffices to
show the two inequalities,
\begin{align}
 \label{e.osc.ineq1}
 \Zbr(t)(A)   & \le\Zbr(t+h)(A^\varepsilon)+ \varepsilon, \\
 \label{e.osc.ineq2}
 \Zbr(t+h)(A) & \le\Zbr(t)(A^\varepsilon)+\varepsilon.
\end{align}

To show \eqref{e.osc.ineq1}, use the definition of $\Omega^r_1$ to write
\begin{equation}
\label{e.osc.ineq1pfa}
\begin{split}
    \Zbr(t)(A) & \le \Zbr(t)(L_\kappa)
      + \Zbr(t)\left(A\cap L_\kappa^c\right)\\
     &\le \frac{\varepsilon}{4}+ \Zbr(t)\left(A\cap L_\kappa^c\right).
\end{split}
\end{equation}
Let $I=\{s\in [t,t+h]:\langle 1,\Zbr(s)\rangle<\varepsilon/4\}$.  Suppose
$I=\varnothing$. Then $\langle 1,\Zbr(s)\rangle\ge\varepsilon/4$ for all
$s\in[t,t+h]$, which implies that
\begin{equation}\label{e.osc.case2a}
  \left\|(\Sr,h)\right\|\le \int_t^{t+\delta}\langle 1,\Zbr(s)\rangle^{-1}\,{\rm d}s
    +\delta
    \le \frac{4\delta}{\varepsilon} +\delta<\min\{\varepsilon,\kappa\}.
\end{equation}
Consequently, $(x,y)\in A\cap L_\kappa^c$ implies $(x,y)-(\Sr,h)\in
A^\varepsilon$, and so
\begin{equation}
    A\cap L_\kappa^c\subset
    A^\varepsilon+(\Sr,h).
\label{e.osc.contain}
\end{equation}
Deduce from \eqref{e.osc.ineq1pfa} that
\begin{equation}
\nonumber
\Zbr(t)(A)  \le \frac{\varepsilon}{4}+ \Zbr(t)\left(A^\varepsilon+(\Sr,h)\right).
\end{equation}
Apply \eqref{e.dynm.bound2} to get
\begin{equation}
\label{e.osc.I01}
\Zbr(t)(A)  \le \frac{\varepsilon}{4}+ \Zbr(t+h)\left(A^\varepsilon\right).
\end{equation}
Suppose $I\neq\varnothing$ and let $\tau =\inf I$. Then $\langle
1,\Zbr(\tau)\rangle\le\varepsilon/4$ by right continuity.  Since $\langle
1,\Zbr(s)\rangle\ge\varepsilon/4$ for all $s\in [t,\tau)$,
\begin{equation}\label{e.osc.Stau}
 \left\|(\overline{S}^{r} (t,\tau),\tau-t)\right\|
	\le\int_t^\tau\langle 1,\Zbr(s)\rangle^{-1}\,{\rm d}s+\delta
  \le  \f{4\delta}{\varepsilon}+\delta
    <\kappa.
\end{equation}
By \eqref{e.osc.ineq1pfa} and \eqref{e.osc.Stau},
\begin{equation}
\nonumber
 \Zbr(t)(A) \le \frac{\varepsilon}{4}+\Zbr(t)(L_\kappa^c)
 	 \le \frac{\varepsilon}{4}+\Zbr(t)\left(
         \overline{\BR}_+^2+(\overline{S}^{r}(t,\tau),\tau-t)\right).
\end{equation}
Apply \eqref{e.dynm.bound2} to get
\begin{equation}
\label{e.osc.In01}
  \Zbr(t)(A) \le \frac{\varepsilon}{4}+\Zbr(\tau)\left( \overline \BR_+^2 \right)
     \le \frac{\varepsilon}{2}.
\end{equation}
So \eqref{e.osc.ineq1} follows because either \eqref{e.osc.I01}
or \eqref{e.osc.In01} holds.

To show \eqref{e.osc.ineq2}, use \eqref{e.dynm.bound1} and the definitions
of $\Omega^r_2$ and $\delta$ to obtain
\begin{align}
\label{e.osc.ineq2pfa}
\begin{split}
  \Zbr(t+h)(A) & \le \Zbr(t)\left(A+(\overline S^r(t,t+h),h)\right)
		+ \Lbr(t,t+h)(\overline \BR_+^2) \\
  & \le \Zbr(t)\left(A+(\overline S^r(t,t+h),h)\right) +\frac{\varepsilon}{4}.
\end{split}
\end{align}
If $I=\varnothing$, then \eqref{e.osc.case2a} implies that $A+(\overline
S^r(t,t+h),h)\subset A^\varepsilon$. So \eqref{e.osc.ineq2pfa} yields
\[%\label{e.osc.}
\Zbr(t+h)(A) \le \Zbr(t)\left(A^\varepsilon\right) +\frac{\varepsilon}{4}.
\]
If $I\ne\varnothing$, then by \eqref{e.dynm.bound1}, the definition of
$\Omega^r_2$ and the choice of $\delta$,
\[%\label{e.osc.case1a}
  \Zbr(t+h)(A)\le  \Zbr(\tau)(\overline \BR_+^2)+\Lbr (\tau,t+h)(\overline \BR_+^2)
    \le \frac{\varepsilon}{4} + 2\lambda^*\delta\le \frac{\varepsilon}{2}.
\]
In both cases, \eqref{e.osc.ineq2} holds.
Conclude from \eqref{e.osc.ineq1} and
\eqref{e.osc.ineq2} that
\begin{equation}
\nonumber
 \dw{\Zbr(t),\Zbr(t+h)}\le \varepsilon.
\end{equation}
Since $t\in[0,T-\delta]$ and $h\in[0,\delta]$ were arbitrary,
\begin{equation}
\nonumber
\ww{\Zbr\cd,\delta}\le\varepsilon,
\end{equation}
which implies that $\omega\in\Omega^r_*$.
\end{proof}

\section{Limiting Fluid Equations}\label{s.limits}
\setcounter{equation}{0}

This section contains the proof of Theorem \ref{t.main}. Tightness
of the sequence $\{\Zbr\cd\}$ follows immediately from Lemmas
\ref{l.compact} and \ref{l.osc}.  Since $\{\Zbr\cd\}$ is tight, there
exists a subsequence $\{q\}\subset\R$ and a process $\Z\cd$ in $\DM$
such that $\overline{\Z}^q\cd\Ar\Z\cd$ as $q\ar\infty$.  We must show that
$\Z\cd$ is almost surely a measure valued fluid model solution for the
data $(\lambda,\vartheta,\zeta_0)$. This is accomplished by Lemmas
\ref{l.lfp-positive} and \ref{l.fl-props}, and Theorem \ref{t.fl}
below.  Finally, if \eqref{a.lip} holds, then a measure valued fluid
model solution for $(\lambda,\vartheta,\zeta_0)$ is unique by Theorem
\ref{t.lip}. In this case, the law of the limit point $\Z\cd$ is unique
and so $\Zbr\cd\Ar\Z\cd$ as $r\ar\infty$.

Let $Z\cd=\langle 1,\Z\cd\rangle$ be the total mass process for $\Z\cd$,
and let $S(u,v)=\int_u^v\frac{1}{Z(s)}\,{\rm d}s$ for all $v\ge u\ge 0$.
To show that $\Z\cd$ is almost surely a measure valued fluid model
solution, note first that $\Z\cd$ is almost surely continuous by Lemma
\ref{l.osc}. Note also that, by \eqref{a.initconv}, $\Z(0)=\zeta_0$
almost surely. It remains to show that properties \eqref{prop-a} and
\eqref{prop-b} of Definition \ref{d.fms} are satisfied almost surely by
$\Z\cd$. The next result establishes \eqref{prop-a}.

\begin{lem}\label{l.lfp-positive}
Almost surely, for all $a>0$,
\begin{equation}
	\inf_{t>a}Z(t)>0.
\label{e.lfp-positive}
\end{equation}
\end{lem}
\begin{proof}
Take $t>0$. Pick a constant $a<t$ small enough such that the marginal
distribution of $D$ is continuous at $a$, take $m<\infty$ such that
the marginal distribution of $B$ is continuous at $m$, and such that
$\lambda\expect{B1_{\{D>a, B<m\}}}>1$.  By dominated convergence,
\[
\lim_{r\rightarrow\infty} 
   \lambda^r\expect{B_1^r 1_{\{D_1^rr^{-1} > a ; 
       B_1^r\leq m \}}}=\lambda \expect{B1_{\{D>a;B\leq m\}}}.
\]
%Using the dynamic equation we get the inequality
%\begin{equation}
%\nonumber
%  Z^r(rt) \geq \sum_{i=E^r(r(t-a))+1}^{E^r(rt)}
%      1^+_{\overline \BR^2_+} \left(B^r_i(rt), D^r_i(rt) \right).
%\end{equation}
%Now replace $B_i^r$
%($i\in \{E^r(r(t-a))+1,\ldots,E^r(rt)\}$) by $B_i^r 1_{\{D_i^r>a; B_i^r\leq m\}}=:B_{i,a,m}^r$.
%This yields
%\begin{equation}
%\nonumber
%  Z^r(rt) \geq
%    \sum_{i=E^r(r(t-a))+1}^{E^r(rt)}
%      1^+_{\overline \BR^2_+} \left(B^r_i(t) \right).
%\end{equation}
Compare the original system with an ordinary PS queues having arrival rate
$\lambda^r_{a,m}= \lambda^r \pr{D^r>ra;B^r<m}$ and service times $B_{i,a,m}^r$,
which are distributed as $B_i^r \mid D_i^r>ra; B_i^r\leq m$.  Suppose
that this PS queue is empty at time $r(t-a)$, and let $\acute{Z}^r(t)$
be the queue length in this PS queue at time $ra$.

Observe that the number of arrivals in the modified PS queue between
time $r(t-a)$ and time $rt$ is less than or equal to the number of
arrivals in that interval in the original PS queue with impatience.
Furthermore, if one of the jobs that arrived in the original PS
queue after time $r(t-a)$ departs before time $rt$, then this must
also be the case in the modified PS queue, since that PS queue had a
service rate which was at least as large as in the original PS queue.
These considerations imply that $Z^r(rt)\geq \acute{Z}^r(rt)$.  Since the
modified queue is still overloaded, and no customer departed because of
impatience, and the modified arrival process is still a renewal process,
the evolution of the modified system between time $r(t-a)$ and $rt$
has the same law as that of an overloaded $GI/GI/1$ PS queue starting
at $0$, in the time interval $[0,ra]$.

Since the service times in our modified system are bounded, the
means converge. The assumptions in \cite{PuStoWi2004} are therefore
valid, and it follows that there exists a constant $k_a>0$ such
that $\lim_{r\rightarrow\infty} \acute{Z}^r(rt)/r=k_a$ almost surely.
Consequently, we have $\liminf_{r\rightarrow\infty} \overline Z^r(t) \geq m_a$
almost surely, which implies the assertion. \end{proof}

Before establishing property \eqref{prop-b} of Definition \ref{d.fms}, 
the following result is needed. 

\begin{lem}
\label{l.fl-props}
Almost surely, for all $C\in\mathcal{C}$ and $t\ge0$,
\begin{equation}
\label{e.lfp-noatoms}
\Z(t)(\partial_C) =0.
\end{equation}
\end{lem}

\begin{proof}   Let $T>0$. It suffices to show the statement for all
$t\in[0,T]$.  Let $\{\eta_n\}\subset(0,1)$ be a sequence such that
$\sum_{n=1}^\infty\eta_n<\infty$. By Lemma \ref{l.noatoms}, there exists
a null sequence of positive reals $\{\kappa_n\}$ such that, for each
fixed $n$,
\begin{equation}
    \liminf_{q\ar\infty}\BP^q\left(
    \sup_{t\in[0,T]}\sup_{C\in\mathcal{C}}
    \overline{\Z}^q(t)(\partial_C^{\kappa_n})\le\frac{1}{n} \right)
    \ge 1-\eta_n.
\label{e.fl-props.kappan}
\end{equation}
For each $n\in\BN$, let $\bM_n=\{\xi\in\bM: \sup_{C\in\mathcal{C}}
\xi(\partial_C^{\kappa_n}) \le{1}/{n}\}$. If a sequence
$\{\xi_i\}\subset\bM_n$ converges weakly to $\xi$, then for each open set
$\partial_C^{\kappa_n}$, the Portmanteau theorem yields
\begin{equation}
    \xi(\partial_C^{\kappa_n})\le
    \limsup_{i\ar\infty}\xi_i(\partial_A^{\kappa_n})
    \le \frac{1}{n}.
\nonumber
\end{equation}
Thus, $\xi\in\bM_n$ and $\bM_n$ is closed. By
definition of the Skorohod $J_1$-topology, the set
$\mathbf{D}_n^T=\{\zeta\cd\in\DM:\zeta(t)\in\bM_n\text{ for all
$t\in[0,T]$}\}$ is also closed. Apply the Portmanteau theorem and
\eqref{e.fl-props.kappan} to obtain
\begin{equation}
    \BP\left(\Z\cd\in\mathbf{D}_n^T\right)
    \ge \liminf_{q\ar\infty}\BP^q\left( \overline{\Z}^q\cd\in\mathbf{D}_n^T\right)
     \ge 1-\eta_n.
\nonumber
\end{equation}
By the Borel-Cantelli lemma,
\begin{equation}
    \BP\left( \bigcup_{k=1}^\infty\bigcap_{n=k}^\infty
    \left\{\Z\cd\in \mathbf{D}_n^T\right\} \right)=1.
\nonumber
\end{equation}
Thus, there exists a finite random variable $N$ such that, almost surely,
\begin{equation}
\label{small.mass.bert}
    \sup_{t\in[0,T]}\sup_{C\in\mathcal{C}}\Z(t)(\partial_C^{\kappa_n})
    \le \frac{1}{n}, \quad \text{for all $n>N$}.
\end{equation}
Since $\partial_C\subset \partial_C^{\kappa_n}$ for all $C\in\mathcal{C}$ and
$n\in\BN$, conclude that almost surely,
\begin{equation}
\notag
    \sup_{t\in[0,T]}\sup_{C\in\mathcal{C}}\Z(t)(\partial_C) =0.
\end{equation}
\end{proof}

We now establish property \eqref{prop-b}.  Recall that $Z(t)=\langle
1,\Z(t)\rangle$ for all $t\ge0$, and $S(u,v)=\int_u^v{1}/{Z(s)}\,{\rm d}s$
for all $v\ge u\ge0$.

\begin{thm}\label{t.fl} Almost surely, the process $\Z\cd$ satisfies
\begin{equation}
\label{e.fl.fmeq}
     \Z(t)\left(A\right)  =\Z(0)\left(A+(S(0,t),t)\right)
       +\lambda\int_0^t \vartheta\left(A+(S(s,t),t-s)\right)\,{\rm d}s,
\end{equation}
for all $t\ge 0$ and $A\in\mathscr{B}$.
\end{thm}

\begin{proof} Let $T>0$. It suffices to show that almost surely, \eqref{e.fl.fmeq}
holds for all $t\in[0,T]$ and all $A\in\mathscr{B}$. For each $r\in\R$,
define the random variable
\begin{equation}
 X^r_T=\sup_{C\in\mathcal{C}}\sup_{0\le s\le t\le T}
  \left|\Lbr(s,t)(C)-\lambda^r(t-s)\breve{\vartheta}^r(C)\right|.
\label{e.fl.Xdef}
\end{equation}
By Lemma \ref{l.glivenko}, $X_T^q\overset{\BP^q}{\lar}0$ as $q\ar\infty$.
Since the limit is deterministic, this convergence is joint with the
convergence $\overline{\Z}^q\cd\Ar\Z\cd$. Using the Skorohod representation
theorem, assume without loss of generality that $\{\overline{\Z}^q\cd,X_T^q\}$
and $\Z\cd$ are defined on a common probability space such that 
\begin{equation}
\label{e.fl.skorohod}
(\overline{\Z}^q\cd,X_T^q)\ar(\Z\cd,0),\quad\text{almost surely.}
\end{equation}
The conclusions of Lemmas \ref{l.lfp-positive} and \ref{l.fl-props}
hold almost surely as well. Assume for the remainder of the proof
that all random objects are evaluated on the event of probability one
such that $\Z\cd$ is continuous, and such that \eqref{e.lfp-positive},
\eqref{e.lfp-noatoms} and \eqref{e.fl.skorohod} hold.

Fix $t\in[0,T]$ and $C\in\mathcal{C}$.  An extension to all Borel sets
$A\in\mathscr{B}$ will be made at the end.  For each $q$, \eqref{e.dyn}
yields
\begin{equation}
\label{e.fl.ZB-dyn}
  \overline{\Z}^q(t)(C)=\overline{\Z}^q(0)
           \left(C+(\overline{S}^{q}(0,t),t)\right)
     + \frac{1}{q}\sum_{i=1}^{q\overline{E}^{q}(t)}
    1^+_C\left(\overline{B}^q_i(t),\overline{D}^q_i(t)\right).
\end{equation}
We will obtain \eqref{e.fl.fmeq} from \eqref{e.fl.ZB-dyn} by
letting $q\ar\infty$.  The convergence in the first component
of \eqref{e.fl.skorohod} is in the Skorohod $J_1$-topology on
$\mathbf{D}([0,\infty),\M)$. However, since $\Z\cd$ is continuous,
\begin{equation}\label{e.fl.wk}
  \overline{\Z}^q(s)\wk\Z(s), 
	\quad \text{for all $s\in[0,t]$.}
\end{equation}
Since $\overline{Z}^q\cd=\langle 1,\overline{\Z}^q\cd\rangle$ and $Z\cd=\langle
1,\Z\cd\rangle$, this implies that
\begin{equation}
\label{e.fl.wk1}
\lim_{q\ar\infty}\left\|\overline{Z}^q\cd -Z\cd\right\|_t= 0.
\end{equation}
For all $t\ge v\ge u>0$, \eqref{e.lfp-positive} implies that
$\inf_{s\in[u,v]}Z(s)>0$, and so the bounded convergence theorem yields
\begin{equation}
\label{e.fl.convergeS}
\begin{split}
\lim_{q\ar\infty}\overline{S}^q(u,v)
 & = \lim_{q\ar\infty}\int_u^v
       \frac{1}{\overline{Z}^q(s)}\,{\rm d}s \\
 & = \int_u^v\frac{1}{Z(s)}\,{\rm d}s \\
 & = S(u,v).
\end{split}
\end{equation}
If $Z(0)\neq0$, then \eqref{e.fl.convergeS} holds for $u=0$ as well,
because then $\inf_{s\in[0,v]}Z(s)>0$.  If $Z(0)=0$, then $S(0,v)=\infty$
and $\overline{S}^q(0,v)\ar\infty$ as $q\ar\infty$.

Suppose that $Z(0)\neq0$ and let $\varepsilon>0$. By
\eqref{e.fl.convergeS}, there exists a $q_\varepsilon\in\R$ such that
$\overline S^q(0,t)\in ((\overline{S}(0,t)-\varepsilon)^+, \overline S(0,t)+\varepsilon)$
for $q>q_\varepsilon$.  Deduce from the shape of the set $C$,
\eqref{e.fl.wk}, and \eqref{e.lfp-noatoms} that
\begin{align*}
\limsup_{q\rightarrow\infty} 
     \overline{\Z}^q(0)\left(C+(\overline{S}^q(0,t),t)\right) &\le
\overline{\Z}(0)\left(C+((\overline S(0,t)-\varepsilon)^+,t)\right),\\
\liminf_{q\rightarrow\infty} 
     \overline{\Z}^q(0)\left(C+(\overline{S}^q(0,t),t)\right) &\ge
\overline{\Z}(0)\left(C+(\overline S(0,t)+\varepsilon,t))\right).
\end{align*}
By $\eqref{e.lfp-noatoms}$, letting $\varepsilon\ar0$ yields
\begin{equation}
\label{e.fl.00}
\lim_{q\ar\infty}\overline{\Z}^q(0)\left(C+(\overline{S}^q(0,t),t)\right)
= \overline{\Z}(0)\left(C+(\overline S(0,t),t)\right).
\end{equation}
If $Z(0)=0$, then \eqref{e.fl.00} holds trivially because the left side
is bounded above by $\lim_{q\ar\infty}\langle 1,\overline{\Z}^q(0)\rangle=0$ by
\eqref{e.fl.wk1}. Combining with \eqref{e.fl.wk} and \eqref{e.lfp-noatoms}
for $\overline{\Z}^q(t)$, implies that, as $q\ar\infty$,
\[%\label{e.lfp.1}
\overline{\Z}^q(t)(C)-\overline{\Z}^q(0)\left(C+(\overline{S}^q(0,t),t)\right)
\ar\Z(t)(C)-\Z(0)\left(C+(S(0,t),t)\right).
\]

Let $I^q$ denote the second right hand term in \eqref{e.fl.ZB-dyn}.
Let $\delta>0$ and let $\eta\in(0,t)$.  Since $\overline{S}^q(s,t)$
is decreasing in $s$ and $S(\cdot,t)$ is continuous on $[\eta,t]$,
\eqref{e.fl.convergeS} implies that $\overline{S}^q(\cdot,t)\ar S(\cdot,t)$
uniformly on $[\eta,t]$.  That is, there exists $q_\delta\in\R$
such that
\begin{equation}
\label{e.fl.uniformS}
 \sup_{s\in[\eta,t]}\left|\overline{S}^q(s,t)
   - S(s,t)\right|\le\delta,
   \qq \text{for all $q>q_\delta$}.
\end{equation}
Let $D_\vartheta(\mathcal{C}) =
\{C\in\mathcal{C}:\vartheta(\partial_C)\ne0
\}$. Note that $D_\vartheta(\mathcal{C})$ is
countable because $\vartheta(\cdot\times\overline{\BR}_+)$ and
$\vartheta(\overline{\BR}_+\times\cdot)$ are probability measures.  Since
$Z(u)>0$ for all $u\in[\eta,t]$, the function $S(s,t)$ is strictly
decreasing in $s$ on $[\eta,t]$. Thus,
\begin{equation}
\nonumber
    D_\vartheta(S)=\left\{s\in [\eta,t]: C+(S(s,t)\pm 2\delta,t-s) \in
    D_\vartheta(\mathcal{C})\right\}
\end{equation}
is also countable.  For each integer $N>1$, let
$\eta=t^N_0<t^N_1<\dots<t^N_N=t$ be a partition of $[\eta,t]$ such
that $t^N_j\notin D_\vartheta(S)$ for all $j=1,\ldots,N-1$, and such
that $\max_{j\le N-1}(t^N_{j+1}-t^N_j)\ar 0$ as $N\ar\infty$.  Then
\[%\label{e.fl.Iq}
    I^q = \frac{1}{q}\sum_{i=1}^{q\overline{E}^q(\eta)}
    1^+_C\left(\overline{B}^q_i(t),\overline{D}^q_i(t)\right) 
   +   \sum_{j=0}^{N-1}
    \frac{1}{q}\sum_{i=q\overline{E}^q(t^N_j)+1}^{q\overline{E}^q(t^N_{j+1})}
    1^+_C\left(\overline{B}^q_i(t),\overline{D}^q_i(t)\right).
\]
Note that the first right hand term is bounded above by
$\overline{\mathcal{L}}^q(0,\eta)(\overline{\BR}_+^2)$.  Suppose that $t^N_j\le
U^q_iq^{-1}\le t^N_{j+1}$, for some $q>q_\delta$, some $j\le N-1$,
and some $i\in\{q\overline{E}^q(\eta)+1,\ldots,q\overline{E}^q(t)\}$.  Then  by
\eqref{e.fl.uniformS},
\begin{equation}
\label{e.fl.sandwich}
     S(t^N_{j+1},t)-\delta \le
    \overline{S}^q(U^q_iq^{-1},t)
    \le S(t^N_j,t)+\delta.
\end{equation}
By definition, 
\begin{equation}
\notag
 \left(\overline{B}^q_i(t),\overline{D}^q_i(t)\right)
  =\left(B^q_i-\overline{S}^q(U^q_iq^{-1},t),
   D^q_iq^{-1}-(  t-U^q_iq^{-1})\right).
\end{equation}
So for $q>q_\delta$, \eqref{e.fl.sandwich} and the inequalities
$1_C(\cdot-\delta,\cdot)\le 1^+_C(\cdot,\cdot)\le 1_C(\cdot+\delta,\cdot)$
yield
\begin{align}
\notag
\begin{split}
     1^+_C\left(\overline{B}^q_i(t),\overline{D}^q_i(t)\right) & \ge
        1_C\left(B^q_i-(S(t^N_j,t)+2\delta),
                D^q_iq^{-1}-(t-t^N_j)\right); \\
     1^+_C\left(\overline{B}^q_i(t),\overline{D}^q_i(t)\right) &\le
        1_C\left(B^q_i-(S(t^N_{j+1},t)-2\delta),
            D^q_iq^{-1}-(t-t^N_{j+1})\right).
\end{split}
\end{align}
This yields, for $q>q_\delta$,
\begin{equation}
\nonumber
%\label{e.fl.pairBounds}
\begin{split}
 I^q & \ge \sum_{j=0}^{N-1}
    \frac{1}{q}\sum_{i=q\overline{E}^q(t^N_j)+1}^{q\overline{E}^q(t^N_{j+1})}
 1_C\left(B^q_i-(S(t^N_j,t)+2\delta),
    D^q_iq^{-1}-(t-t^N_j)\right); \\
 I^q & \le  \overline{\mathcal{L}}^q(0,\eta)(\overline{\BR}_+^2)\\
	&\quad
   	 +\sum_{j=0}^{N-1}
    \frac{1}{q}\sum_{i=q\overline{E}^q(t^N_j)+1}^{q\overline{E}^q(t^N_{j+1})}
   1_C\left(B^q_i-(S(t^N_{j+1},t)-2\delta),
      D^q_iq^{-1}-(t-t^N_{j+1})\right).
\end{split}
\end{equation}
Rewrite as
\begin{equation}
\label{e.lfp.bounds}
\begin{split}
  I^q  &\ge \sum_{j=0}^{N-1}
    \overline{\mathcal{L}}^q(t^N_j,t^N_{j+1})\left(C+\left(S(t^N_j,t)+2\delta,
       t-t^N_j\right)\right); \\
      I^q &\le \overline{\mathcal{L}}^q(0,\eta)(\overline{\BR}_+^2)
      +\sum_{j=0}^{N-1}
    \overline{\mathcal{L}}^q(t^N_j,t^N_{j+1})\left(
       C+\left(S(t^N_{j+1},t)-2\delta,t-t^N_{j+1}\right)\right).
\end{split}
\end{equation}
By \eqref{e.fl.Xdef} and \eqref{e.lfp.bounds}, $q>q_\delta$ implies that
\begin{align*}
  I^q & \ge \sum_{j=0}^{N-1}
    \left(\lambda^q(t^N_{j+1}-t^N_j)
    \breve{\vartheta}^q\left(C+\left(S(t^N_j,t)+2\delta,t-t^N_j\right)\right)
     -X_T^q\right); \\
      I^q & \le \lambda^q\eta+X_T^q+
   \sum_{j=0}^{N-1}
    \left(\lambda^q(t^N_{j+1}-t^N_j)
    \breve{\vartheta}^q\left(C+\left(S(t^N_{j+1},t)-2\delta,t-t^N_{j+1}\right)\right)
     +X_T^q\right).
\end{align*}
By \eqref{e.fl.skorohod}, and since $t^N_j\not\in D_\vartheta(S)$ for all
$j=1,\ldots,N-1$,
\begin{equation}\label{e.fl.Bounds}
\begin{split}
  \liminf_{q\ar\infty}I^q
  & \ge  \lambda\sum_{j=0}^{N-1}
    (t^N_{j+1}-t^N_j)
    \vartheta\left(C+\left(S(t^N_j,t)+2\delta,t-t^N_j\right)\right); \\
  \limsup_{q\ar\infty}I^q
   & \le \lambda\eta + \lambda\sum_{j=0}^{N-1}
    (t^N_{j+1}-t^N_j)
    \vartheta\left(C+\left(S(t^N_{j+1},t)-2\delta,t-t^N_{j+1}\right)\right).
\end{split}
\end{equation}
For $s\in[\eta,t]$ such that $s\notin D_\vartheta(S)$ the bounded
convergence theorem implies that
\begin{equation}
\label{e.fl.bc}
\begin{split}
  &\lim_{N\to+\infty} \sum_{j=0}^{N-1}1_{[t^N_j,t^N_{j+1})}(s)
    \vartheta\left(C+\left(S(t^N_j,t)+2\delta,t-t^N_j\right)\right)
    \\& \qquad\qquad\qquad=\vartheta\left(C+\left(S(s,t)+2\delta,t-s\right)\right);\\
  &\lim_{N\to+\infty}  \sum_{j=0}^{N-1}1_{[t^N_j,t^N_{j+1})}(s)
    \vartheta\left(C+\left(S(t^N_{j+1},t)-2\delta,t-t^N_{j+1}\right)\right)
    \\& \qquad\qquad\qquad  =\vartheta\left(C+\left(S(s,t)-2\delta,t-s\right)\right).
\end{split}
\end{equation}
Thus, the convergence in \eqref{e.fl.bc} holds for almost every
$s\in[\eta,t)$. Let $N\ar\infty$ in \eqref{e.fl.Bounds} and conclude from
\eqref{e.fl.bc} and the bounded convergence theorem that
\begin{equation}
\label{e.lfp.lastBounds}
\begin{split}
  \liminf_{q\ar\infty}I^q
  & \ge  \lambda\int_\eta^t
    \vartheta\left(C+\left(S(s,t)+2\delta,t-s\right)\right)\,{\rm d}s; \\
  \limsup_{q\ar\infty}I^q
   & \le \lambda\eta+\lambda\int_\eta^t
    \vartheta\left(C+\left(S(s,t)-2\delta,t-s\right)\right)\,{\rm d}s.
\end{split}
\end{equation}
Let $\delta\ar0$ in \eqref{e.lfp.lastBounds}. Since
$D_\vartheta(\mathcal{C})$ is countable, both integrands in
\eqref{e.lfp.lastBounds} converge almost everywhere on $[\eta,t]$ to
$\vartheta\left(C+\left(S(s,t),t-s\right)\right)$.  Thus,
\begin{equation}
\nonumber
\begin{split}
  \liminf_{q\ar\infty}I^q
  & \ge  \lambda\int_\eta^t
    \vartheta\left(C+\left(S(s,t),t-s\right)\right)\,{\rm d}s; \\
  \limsup_{q\ar\infty}I^q
   & \le \lambda\eta+\lambda\int_\eta^t
    \vartheta\left(C+\left(S(s,t),t-s\right)\right)\,{\rm d}s.
\end{split}
\end{equation}
Let $\eta\ar0$ to conclude that
\begin{equation}
\notag
    \lim_{q\ar\infty} I^q = \lambda\int_0^t
     \vartheta\left(C+\left(S(s,t),t-s\right)\right)\,{\rm d}s.
\end{equation}

This proves \eqref{e.fl.fmeq} for all $t\in[0,T]$ and $C\in\mathcal{C}$.  To extend to all
$A\in\mathscr{B}$,  let   $\mathcal{C}'$  be  the  set  of   $A\in\mathscr{B}$  for  which
\eqref{e.fl.fmeq}   holds.    Observe   that   $\mathcal{C}'$   is   a   $\lambda$-system:
$\overline{\BR}_+^2\in\mathcal{C}'$       because       $\overline{\BR}_+^2\in\mathcal{C}$;       if
$\{A_n\}\subset\mathcal{C}'$  satisfies  $A_n  \uparrow  A$, then  $A\in\mathcal{C}'$;  if
$A_1\subset   A_2$  are   elements   of  $\mathcal{C}'$,   then   $A_2\setminus  A_1   \in
\mathcal{C}'$.    Observe    also   that    $\mathcal{C}$    is    a   $\pi$-system:    if
$C_1,C_2\in\mathcal{C}$,       then        $C_1\cap       C_2\in\mathcal{C}$.        Since
$\mathcal{C}\subset\mathcal{C}'$  and the $\sigma$-algebra  generated by  $\mathcal{C}$ is
equal  to  $\mathscr{B}$,  it   follows  that  $\mathcal{C}'=\mathscr{B}$  by  the  Dynkin
$\pi\lambda$-theorem (see for example \cp{Bi1986}).  \end{proof}

\providecommand{\bysame}{\leavevmode\hbox to3em{\hrulefill}\thinspace}
\providecommand{\MR}{\relax\ifhmode\unskip\space\fi MR }
\providecommand{\MRhref}[2]{%
  \href{http://www.ams.org/mathscinet-getitem?mr=#1}{#2}
}
\providecommand{\href}[2]{#2}


\begin{thebibliography}{10}

\bibitem{B57}
D.~Barrer, \emph{Queueing with impatient customers and ordered service},
  Operations Research \textbf{5} (1957), 650--656.

\bibitem{Bi1986}
Patrick Billingsley, \emph{Probability and measure}, 2 ed., John Wiley \& Sons,
  Inc.\, New York, 1986.

\bibitem{BonMas2001}
Thomas Bonald and Laurent Massouli\'{e}, \emph{Impact of fairness on {I}nternet
  performance}, Proceedings of ACM Sigmetrics 2001, 2001, pp.~82--91.

\bibitem{BR03}
Thomas Bonald and James Roberts, \emph{Congestion at flow level and the impact
  of user behaviour}, Computer Networks \textbf{42} (2003), 521--536.

\bibitem{BT98}
N.~Boots and Tijms H., \emph{A multi-server queueuing system with impatient
  customers}, Management Science \textbf{45} (1999), 444--448.

\bibitem{Bra2005}
M.~Bramson, \emph{Stability of networks for max-min fair routing}, Presentation
  at the 13th INFORMS Applied Probability Conference, Ottawa, 2005.

\bibitem{CPRW94}
E.~Coffman, A.~Puhalskii, M.~Reiman, and P.~Wright, \emph{Processor shared
  buffers with reneging}, Performance Evaluation \textbf{19} (1994), 25--46.

\bibitem{deVKonLee2001}
Gustavo de~Veciana, Takis Konstantopoulos, and Tae-Jin Lee, \emph{Stability and
  performance analysis of networks supporting elastic services}, IEEE/ACM
  Trans. Netw. \textbf{9} (2001), no.~1, 2--14.

\bibitem{DoLeSh2001}
Bogdan Doytchinov, John Lehoczky, and Steven Shreve, \emph{Real-time queues in
  heavy traffic with earliest-deadline-first queue discipline}, Annals of
  Applied Probability \textbf{11} (2001), no.~2, 332--378.

\bibitem{GKM02}
N.~Gans, G.~Koole, and A.~Mandelbaum, \emph{Telephone call centers: Tutorial,
  review, and research prospects}, Manufacturing \& Service Operations
  Management \textbf{5} (2002), 79--141.

\bibitem{Gromoll:02}
Christian Gromoll, Philippe Robert, Bert Zwart, and Richard Bakker, \emph{The
  impact of reneging in processor sharing queues}, ACM-Sigmetrics (Saint Malo),
  ACM/IFIP WG 7.3, June 2006.

\bibitem{GrKru2004}
H.~C. Gromoll and \L. Kruk, \emph{Heavy traffic analysis of a real-time
  processor sharing queue}, to appear, 2006.

\bibitem{GrWi2005}
H.~C. Gromoll and R.~J. Williams, \emph{Fluid approximation for an
  \uppercase{I}nternet congestion control model with fair bandwidth sharing and
  general document size distributions}, Preprint, 2006.

\bibitem{GRZ04}
Fabrice Guillemin, Philippe Robert, and Bert Zwart, \emph{Tail asymptotics for
  processor-sharing queues}, Advances in Applied Probability \textbf{36}
  (2004), 525--543.

\bibitem{HV93}
J.~Hale and S.~Verduyn~Lunel, \emph{An introduction to functional differential
  equations}, Springer Verlag, New York, 1993.

\bibitem{JeanMarie:02}
Alain Jean-Marie and Philippe Robert, \emph{On the transient behavior of some
  single server queues}, Queueing Systems, Theory and Applications \textbf{17}
  (1994), 129--136.

\bibitem{Ka1976}
Olav Kallenberg, \emph{Random {M}easures}, Academic Press, New York, 1986.

\bibitem{KeWi2004}
F.~P.\ Kelly and R.~J.\ Williams, \emph{Fluid model for a network operating
  under a fair bandwidth sharing policy}, Annals of Applied Probability
  \textbf{14} (2004), 1055--1083.

\bibitem{KMBK04}
P.~Key, L.\ Massouli\'e, A.~Bain, and F.~Kelly, \emph{Fair internet traffic
  integration: Network flow models and analysis}, Annals of Telecommunications
  \textbf{59} (2004), 1338--1352.

\bibitem{KruLehShr2003}
L.~Kruk, J.~Lehoczky, and S.~Shreve, \emph{Second order approximation for the
  customer time in queue distribution under the {FIFO} service discipline},
  Annales UMCS Informatica AI \textbf{1} (2003), 37--48.

\bibitem{KruLehShr2004}
\bysame, \emph{Accuracy of state space collapse for earliest-deadline-first
  queues}, Preprint, 2004.

\bibitem{KruLehShrYeu2001}
L.~Kruk, J.~Lehoczky, S.~Shreve, and S.~Yeung, \emph{Multiple-input
  heavy-traffic real-time queues}, Annals of Applied Probability \textbf{13}
  (2003), no.~1, 54--99.

\bibitem{KruLehShrYeu2004}
\bysame, \emph{Earliest-deadline-first service in heavy-traffic acyclic
  networks}, Annals of Applied Probability \textbf{14} (2004), no.~3,
  1306--1352.

\bibitem{LaBeSr2004}
A.~Lakshmikantha, C.~L. Beck, and R.~Srikant, \emph{Connection level stability
  analysis of the internet using the sum of squares ({SoS}) techniques},
  Conference on Information Sciences and Systems, Princeton, 2004.

\bibitem{Mas2005}
Laurent Massouli\'e, \emph{Structural properties of proportional fairness:
  stability and insensitivity}, Preprint, 2005.

\bibitem{Massoulie}
Laurent Massouli\'e and James Roberts, \emph{Bandwidth sharing: Objectives and
  algorithms}, {INFOCOM }'99. Eighteenth Annual Joint Conference of the {IEEE}
  Computer and Communications Societies, 1999, pp.~1395--1403.

\bibitem{MoWal2000}
J.~Mo and J.~Walrand, \emph{Fair end-to-end window-based congestion control},
  IEEE/ACM Transactions on Networking \textbf{8} (2000), no.~5, 556--567.

\bibitem{PuStoWi2004}
A.~L. Puha, A.~L. Stolyar, and R.~J. Williams, \emph{The fluid limit of an
  overloaded processor sharing queue}, Preprint, 2004.

\bibitem{RoMa2000}
James Roberts and Laurent Massouli\'{e}, \emph{Bandwidth sharing and admission
  control for elastic traffic}, Telecommunication Systems \textbf{15} (2000),
  185--201.

\bibitem{S79}
Robert~E. Stanford, \emph{Reneging phenomena in single channel queues},
  Mathematics of Operations Research \textbf{4} (1979), 162--178.

\bibitem{Stanford}
\bysame, \emph{On queues with impatience}, Advances in Applied Probability
  \textbf{22} (1990), no.~3, 768--769.

\bibitem{VaaWel1996}
Aad van~der Vaart and Jon~A. Wellner, \emph{Weak convergence and empirical
  processes}, Springer-Verlag, New York, 1996.

\bibitem{WG04}
A.~Ward and P.~Glynn, \emph{A diffusion approximation for a markovian queue
  with reneging}, Queueing Systems \textbf{43} (2003), 103--128.

\bibitem{YeuLeh2004}
Shu-Ngai Yeung and John~P. Lehoczky, \emph{Real-time queueing networks in heavy
  traffic with {EDF} and {FIFO} queue discipline}, Preprint, 2004.

\end{thebibliography}
\end{document}